\documentclass[12pt]{article}

\usepackage{amsmath,amssymb,amsthm,amscd,a4wide}

\begin{document}

\newcommand{\End}{{\rm{End}\ts}}
\newcommand{\Hom}{{\rm{Hom}}}
\newcommand{\ch}{{\rm{ch}\ts}}
\newcommand{\non}{\nonumber}
\newcommand{\wt}{\widetilde}
\newcommand{\wh}{\widehat}
\newcommand{\ot}{\otimes}
\newcommand{\la}{\lambda}
\newcommand{\La}{\Lambda}
\newcommand{\De}{\Delta}
\newcommand{\al}{\alpha}
\newcommand{\be}{\beta}
\newcommand{\ga}{\gamma}
\newcommand{\Ga}{\Gamma}
\newcommand{\ka}{\kappa}
\newcommand{\si}{\sigma}
\newcommand{\vp}{\varphi}
\newcommand{\de}{\delta^{}}
\newcommand{\ze}{\zeta}
\newcommand{\om}{\omega}
\newcommand{\hra}{\hookrightarrow}
\newcommand{\ve}{\varepsilon}
\newcommand{\ts}{\,}
\newcommand{\vac}{\mathbf{1}}
\newcommand{\di}{\partial}
\newcommand{\qin}{q^{-1}}
\newcommand{\tss}{\hspace{1pt}}
\newcommand{\Sr}{ {\rm S}}
\newcommand{\U}{ {\rm U}}
\newcommand{\Y}{ {\rm Y}}
\newcommand{\BL}{ {\overline L}}
\newcommand{\BE}{ {\overline E}}
\newcommand{\BP}{ {\overline P}}
\newcommand{\AAb}{\mathbb{A}\tss}
\newcommand{\CC}{\mathbb{C}\tss}
\newcommand{\QQ}{\mathbb{Q}\tss}
\newcommand{\SSb}{\mathbb{S}\tss}
\newcommand{\ZZ}{\mathbb{Z}\tss}
\newcommand{\Z}{{\rm Z}}
\newcommand{\Ac}{\mathcal{A}}
\newcommand{\Lc}{\mathcal{L}}
\newcommand{\Mc}{\mathcal{M}}
\newcommand{\Pc}{\mathcal{P}}
\newcommand{\Qc}{\mathcal{Q}}
\newcommand{\Tc}{\mathcal{T}}
\newcommand{\Sc}{\mathcal{S}}
\newcommand{\Bc}{\mathcal{B}}
\newcommand{\Ec}{\mathcal{E}}
\newcommand{\Hc}{\mathcal{H}}
\newcommand{\Uc}{\mathcal{U}}
\newcommand{\Vc}{\mathcal{V}}
\newcommand{\Wc}{\mathcal{W}}
\newcommand{\Ar}{{\rm A}}
\newcommand{\Ir}{{\rm I}}
\newcommand{\Fr}{{\rm F}}
\newcommand{\Jr}{{\rm J}}
\newcommand{\Rr}{{\rm R}}
\newcommand{\Zr}{{\rm Z}}
\newcommand{\gl}{\mathfrak{gl}}
\newcommand{\middd}{{\rm mid}}
\newcommand{\Pf}{{\rm Pf}}
\newcommand{\Norm}{{\rm Norm\tss}}
\newcommand{\oa}{\mathfrak{o}}
\newcommand{\spa}{\mathfrak{sp}}
\newcommand{\g}{\mathfrak{g}}
\newcommand{\h}{\mathfrak h}
\newcommand{\n}{\mathfrak n}
\newcommand{\z}{\mathfrak{z}}
\newcommand{\Zgot}{\mathfrak{Z}}
\newcommand{\p}{\mathfrak{p}}
\newcommand{\sll}{\mathfrak{sl}}
\newcommand{\agot}{\mathfrak{a}}
\newcommand{\qdet}{ {\rm qdet}\ts}
\newcommand{\Ber}{ {\rm Ber}\ts}
\newcommand{\HC}{ {\mathcal HC}}
\newcommand{\cdet}{ {\rm cdet}}
\newcommand{\tr}{ {\rm tr}}
\newcommand{\str}{ {\rm str}}
\newcommand{\loc}{{\rm loc}}
\newcommand{\Gr}{ {\rm Gr}\tss}
\newcommand{\sgn}{ {\rm sgn}\ts}
\newcommand{\ba}{\bar{a}}
\newcommand{\bb}{\bar{b}}
\newcommand{\bi}{\bar{\imath}}
\newcommand{\bj}{\bar{\jmath}}
\newcommand{\bk}{\bar{k}}
\newcommand{\bl}{\bar{l}}
\newcommand{\Sym}{\mathfrak S}
\newcommand{\fand}{\quad\text{and}\quad}
\newcommand{\Fand}{\qquad\text{and}\qquad}
\newcommand{\vt}{{\tss|\hspace{-1.5pt}|\tss}}
\newcommand{\antiddots}
    {\underset{\displaystyle\cdot\quad\ }
    {\overset{\displaystyle\quad\ \cdot}{\cdot}}}
\newcommand{\dddots}
    {\underset{\displaystyle\quad\ \cdot}
    {\overset{\displaystyle\cdot\quad\ }{\cdot}}}

\renewcommand{\theequation}{\arabic{section}.\arabic{equation}}

\newtheorem{thm}{Theorem}[section]
\newtheorem{lem}[thm]{Lemma}
\newtheorem{prop}[thm]{Proposition}
\newtheorem{cor}[thm]{Corollary}
\newtheorem{conj}[thm]{Conjecture}
\newtheorem*{mthma}{Theorem A}
\newtheorem*{mthmb}{Theorem B}

\theoremstyle{definition}
\newtheorem{defin}[thm]{Definition}

\theoremstyle{remark}
\newtheorem{remark}[thm]{Remark}
\newtheorem{example}[thm]{Example}

\newcommand{\bth}{\begin{thm}}
\renewcommand{\eth}{\end{thm}}
\newcommand{\bpr}{\begin{prop}}
\newcommand{\epr}{\end{prop}}
\newcommand{\ble}{\begin{lem}}
\newcommand{\ele}{\end{lem}}
\newcommand{\bco}{\begin{cor}}
\newcommand{\eco}{\end{cor}}
\newcommand{\bde}{\begin{defin}}
\newcommand{\ede}{\end{defin}}
\newcommand{\bex}{\begin{example}}
\newcommand{\eex}{\end{example}}
\newcommand{\bre}{\begin{remark}}
\newcommand{\ere}{\end{remark}}
\newcommand{\bcj}{\begin{conj}}
\newcommand{\ecj}{\end{conj}}

\newcommand{\bal}{\begin{aligned}}
\newcommand{\eal}{\end{aligned}}
\newcommand{\beq}{\begin{equation}}
\newcommand{\eeq}{\end{equation}}
\newcommand{\ben}{\begin{equation*}}
\newcommand{\een}{\end{equation*}}

\newcommand{\bpf}{\begin{proof}}
\newcommand{\epf}{\end{proof}}

\def\beql#1{\begin{equation}\label{#1}}

\title{\Large\bf Combinatorial bases for
covariant representations\\
of the Lie superalgebra $\gl_{m|n}$}

\author{A. I. Molev}

\date{} 
\maketitle

\vspace{15 mm}

\begin{abstract}
Covariant tensor representations of $\gl_{m|n}$
occur as irreducible components of tensor powers
of the natural $(m+n)$-dimensional representation.
We construct a basis of each covariant representation
and give explicit formulas
for the action of the generators of $\gl_{m|n}$ in
this basis. The basis has the property that
the natural Lie subalgebras $\gl_m$ and $\gl_n$ act by
the classical Gelfand--Tsetlin formulas.
The main role in the construction is played by the fact that
the subspace of $\gl_m$-highest vectors in any
finite-dimensional irreducible representation of $\gl_{m|n}$
carries a structure of an irreducible module
over the Yangian $\Y(\gl_n)$. One consequence is a new proof
of the character formula for the covariant representations
first found by Berele and Regev and by Sergeev.
\end{abstract}


\vspace{35 mm}

\noindent
School of Mathematics and Statistics\newline
University of Sydney,
NSW 2006, Australia\newline
Email: alexander.molev@sydney.edu.au

\newpage

\section{Introduction}\label{sec:int}
\setcounter{equation}{0}

Finite-dimensional
irreducible representations of the Lie
superalgebra $\gl_{m|n}$ over $\CC$
are parameterized by their highest
weights. The highest weight of such a representation
is a tuple $\la$ of complex numbers
of the form $\la=(\la_1,\dots,\la_m\ts|\ts\la_{m+1},\dots,\la_{m+n})$
satisfying the conditions
\beql{findim}
\la_i-\la_{i+1}\in\ZZ_+,\qquad
\text{for}\quad i=1,\dots,m+n-1,\quad i\ne m.
\eeq
We let $L(\la)$ denote the corresponding
representation. It is isomorphic to the unique
irreducible quotient of the
{\it Kac module\/} $K(\la)$ which is defined as a universal
induced module associated with the irreducible
module over the Lie subalgebra $\gl_{m}\oplus \gl_{n}$ with
the highest weight $\la$. Due to the work of Kac~\cite{k:ct, k:rc},
the module $K(\la)$ is irreducible if and only if $\la$ is a
{\it typical\/} weight which thus provides a character formula
for the irreducible representations $L(\la)$ associated
with typical highest weights.

The problem of finding
the characters of the complementary family
of {\it atypical\/} representations remained open
until the work of Serganova~\cite{s:kl}. She produced
character formulas for $L(\la)$ in terms of the characters
of the Kac modules involving
generalized Kazhdan--Lusztig polynomials.
These results were extended and made more explicit
by Brundan~\cite{b:kl}
and Su and Zhang~\cite{s:cd} which led to
character and dimension formulas for all representations
$L(\la)$ and allowed to prove conjectures
on the characters stated by several authors; see
\cite{b:kl} and \cite{s:cd} for more
detailed discussions and references.

In this paper we develop a different approach allowing to employ
the Yangian representation theory to construct representations
of $\gl_{m|n}$ and calculate their characters.
We consider the {\it multiplicity space\/} $L(\la)^+_{\mu}$
spanned by the $\gl_m$-highest vectors in $L(\la)$ of weight $\mu$,
which is isomorphic to
the space of $\gl_m$-homomorphisms $\Hom_{\gl_m}(L'(\mu),L(\la))$, where
$L'(\mu)$ denotes irreducible representation of the Lie algebra
$\gl_m$ with the highest weight $\mu$.
The multiplicity space is a natural module
over the centralizer $\U(\gl_{m|n})^{\gl_m}$. As with
the Lie algebra case, there exists an `almost surjective'
homomorphism from the {\it Yangian\/} $\Y(\gl_n)$
to the centralizer thus yielding an irreducible action of $\Y(\gl_n)$
on $L(\la)^+_{\mu}$. This is analogous to the Olshanski
{\it centralizer construction\/}, originally appeared
in the context of the classical Lie algebras \cite{o:ri, o:ty}
and which has led, in particular, to constructions of weight
bases for representations of classical Lie algebras
\cite{m:yc}.

On the other hand, finite-dimensional irreducible representations
of the Yangian $\Y(\gl_n)$ are classified in terms of their
highest weights or Drinfeld polynomials. Hence, identifying
the parameters of
the representation $L(\la)^+_{\mu}$ of $\Y(\gl_n)$ we
can, in principle, find its character by using the results
on the characters of the Yangian modules;
see \cite[Theorem~15]{a:df} and \cite[Corollary~8.22]{bk:rs}.

We apply the Yangian approach to
a particular family of
{\it covariant tensor representations\/} $L(\la)$ of $\gl_{m|n}$
(which we simply call the {\it covariant representations\/}).
They occur as irreducible components of the tensor powers
of the natural $(m+n)$-dimensional representation $\CC^{m|n}$.
The irreducible decompositions of the tensor powers
and the characters of the covariant representations
were first found by Berele and Regev~\cite{br:hy}
and Sergeev~\cite{s:ta}. The character formulas are expressed
in terms of the {\it supersymmetric Schur polynomials\/}
which admit a combinatorial presentation
in terms of {\it supertableaux\/}.
The covariant representations of $\gl_{m|n}$ include both
typical and atypical representations and share many common
properties with the polynomial representations of $\gl_{m}$
arising in the decomposition of the tensor powers of $\CC^{m}$.

Our main result is an explicit construction of
all covariant representations $L(\la)$
of the Lie superalgebra
$\gl_{m|n}$. We construct
a basis of $L(\la)$ parameterized by supertableaux and
give explicit formulas for the action of the generators of $\gl_{m|n}$.
We do this by using the vector space isomorphism
\beql{decomptp}
L(\la)\cong \bigoplus _{\mu}L'(\mu)\ot L(\la)^+_{\mu},
\eeq
summed over $\gl_m$-highest weights $\mu$.
It turns out that if $L(\la)$ is a covariant representation, then
the $\Y(\gl_n)$-module $L(\la)^+_{\mu}$
is isomorphic to the {\it skew module\/} $\BL(\la')^+_{\mu'}$
arising from the Olshanski centralizer construction applied
to the {\it Lie algebra\/} $\gl_{r+n}$ and its subalgebra $\gl_r$.
The skew modules were introduced and studied by Cherednik~\cite{c:ni}.
Their Drinfeld polynomials
were first calculated by Nazarov
and Tarasov~\cite{nt:ry}, and different proofs were also given
later in \cite{hm:qa} and \cite{m:yt}; see also \cite[Sec. 8.5]{m:yc}.
The skew module $\BL(\la')^+_{\mu'}$ possesses a basis parameterized
by trapezium-like Gelfand--Tsetlin patterns, and hence so does
the multiplicity space $L(\la)^+_{\mu}$. Moreover, the
$\gl_m$-module $L'(\mu)$ admits the classical Gelfand--Tsetlin basis,
so that a basis of $L(\la)$ can be naturally parameterized
by pairs of patterns. An equivalent combinatorial description
of the basis in terms of supertableaux is also given.

To calculate the Yangian highest weight of the representation
$L(\la)^+_{\mu}$ and to
derive the matrix element formulas for the generators of
$\gl_{m|n}$ in this basis we employ a relationship between
the Yangian $\Y(\gl_n)$ and the {\it Mickelsson--Zhelobenko algebra\/}
$\Z(\gl_{m|n},\gl_m)$. Namely, we find
the images of generators of the algebra
$\Y(\gl_n)$ under the composition
\ben
\Y(\gl_n)\to\U(\gl_{m|n})^{\gl_m}\to \Z(\gl_{m|n},\gl_m),
\een
where the second map is a natural homomorphism.
The elements of each of the natural subalgebras $\gl_m$ and $\gl_n$
act in the basis of $L(\la)$ by the classical Gelfand--Tsetlin
formulas, and the action of the odd component of the Lie superalgebra
is also expressed in an explicit form; cf.
\cite[Ch. 9]{m:yc} where a similar approach was used
to construct weight bases for representations
of the classical Lie algebras.

In our arguments we will avoid using the character
formulas of \cite{br:hy} and \cite{s:ta} so we thus obtain a new
proof of those formulas.

Note that analogues of the Gelfand--Tsetlin bases for a certain
class of {\it essentially typical\/}
representations of the Lie superalgebra $\gl_{m|n}$
were given in \cite{p:et} and \cite{tis:gt}. The constructions
rely on the fact that consecutive restrictions of an essentially typical
representation to the subalgebras of the chain
$
\gl_{m|1}\subset \gl_{m|2}\subset\dots\subset\gl_{m|n}
$
are completely reducible thus allowing to apply the
approach of Gelfand and Tsetlin~\cite{gt:fdu} to construct
a basis and to derive explicit formulas for the action
of the generators. In the special case $n=1$
any finite-dimensional irreducible representation
of the Lie superalgebra $\gl_{m|1}$ can be realized via a basis of
Gelfand--Tsetlin type by using the restriction to
the natural Lie subalgebra $\gl_m$; see \cite{p:if}.
In a recent work \cite{sv:gz} these results were extended
to all covariant representations of $\gl_{m|n}$ by constructing
a Gelfand--Tsetlin type basis and providing matrix element
formulas of the generators of $\gl_{m|n}$ in the basis.
That basis appears to have properties quite different
from ours.

\section{Main results}\label{sec:mt}
\setcounter{equation}{0}

Consider the standard basis
$E_{ij}$, $1\leqslant i,j\leqslant m+n$, of the
Lie superalgebra
$\gl_{m|n}$ over $\CC$. The $\ZZ_2$-grading on $\gl_{m|n}$
is defined by setting $\deg E_{ij}=\bi+\bj$, where
we use the notation $\bi=0$ for $1\leqslant i\leqslant m$
and $\bi=1$ for $m+1\leqslant i\leqslant m+n$.
The commutation relations
have the form
\ben
\big[E_{ij},E_{kl}\tss\big]
=\de_{kj}\ts E_{i\tss l}-\de_{i\tss l}\ts
E_{kj}(-1)^{(\bi+\bj)(\bk+\bl)},
\een
where the square brackets denote the super-commutator.

Given a tuple $\la=(\la_1,\dots,\la_m\ts|\ts\la_{m+1},\dots,\la_{m+n})$
of complex numbers, the irreducible representation
$L(\la)$ of $\gl_{m|n}$ is generated
by a nonzero vector $\ze$ (the highest vector)
satisfying the conditions
\begin{alignat}{2}
E_{ij}\ts\ze&=0 \qquad &&\text{for} \quad
1\leqslant i<j\leqslant m+n \qquad \text{and}
\non\\
E_{ii}\ts\ze&=\la_i\ts\ze \qquad &&\text{for}
\quad 1\leqslant i\leqslant m+n.
\non
\end{alignat}
The representation $L(\la)$ is finite-dimensional if and only if
the conditions \eqref{findim} hold. We will suppose
that this is the case as we will only work with
finite-dimensional representations.

The family of the {\it covariant representations\/}
is distinguished by the following conditions:
all components $\la_1,\dots,\la_{m+n}$ of $\la$
are nonnegative integers
and the number $\ell$
of nonzero components among $\la_{m+1},\dots,\la_{m+n}$
does not exceed $\la_m$; see \cite{br:hy}
and \cite{s:ta}. To each highest weight
$\la$ satisfying these conditions, we will associate
the Young diagram $\Ga_{\la}$ containing $\la_1+\dots+\la_{m+n}$ boxes
which is determined by the conditions that
the first $m$ rows of $\Ga_{\la}$ are $\la_1,\dots,\la_m$
while the first $\ell$ columns are
$\la_{m+1}+m,\dots,\la_{m+\ell}+m$. The condition
$\ell\leqslant\la_m$ ensures that $\Ga_{\la}$ is
the diagram of a partition. A {\it supertableau\/} $\La$ of shape
$\Ga_{\la}$ is obtained by filling in the boxes of
the diagram $\Ga_{\la}$ with the numbers $1,\dots,m+n$
in such a way that
\ben
\bal
&\text{the entries weakly increase
from left to right along each row
and down
each column;}\\
&\text{the entries in $\{1,\dots,m\}$ strictly increase
down each column;}\\
&\text{the entries in $\{m+1,\dots,m+n\}$ strictly increase
from left to right along each row.}
\eal
\een

\bex\label{ex:st} The following is a supertableau
of shape $\Ga_{\la}$
associated with the highest weight
$\la=(10,7,4,3\ts|\ts 3,1,0,0,0)$
of $\gl_{4|5}$:

\setlength{\unitlength}{1.5em}
\begin{center}
\begin{picture}(9,8)

\put(0,0){\line(0,1){7}}
\put(1,0){\line(0,1){7}}
\put(2,2){\line(0,1){5}}
\put(3,3){\line(0,1){4}}
\put(4,4){\line(0,1){3}}
\put(5,5){\line(0,1){2}}
\put(6,5){\line(0,1){2}}
\put(7,5){\line(0,1){2}}
\put(8,6){\line(0,1){1}}
\put(9,6){\line(0,1){1}}
\put(10,6){\line(0,1){1}}

\put(0,7){\line(1,0){10}}
\put(0,6){\line(1,0){10}}
\put(0,5){\line(1,0){7}}
\put(0,4){\line(1,0){4}}
\put(0,3){\line(1,0){3}}
\put(0,2){\line(1,0){2}}
\put(0,1){\line(1,0){1}}
\put(0,0){\line(1,0){1}}

\multiput(3,3)(0.4,0){5}{\line(1,0){0.2}}

\multiput(5,3)(0,0.4){6}{\line(0,1){0.2}}

\linethickness{1.2pt}

\put(0,3){\line(1,0){1}}
\put(1,4){\line(1,0){1}}
\put(2,5){\line(1,0){4}}
\put(1,3){\line(0,1){1}}
\put(2,4){\line(0,1){1}}
\put(6,5){\line(0,1){2}}

\put(0.3,6.3){$1$}
\put(1.3,6.3){$1$}
\put(2.3,6.3){$1$}
\put(3.3,6.3){$2$}
\put(4.3,6.3){$2$}
\put(5.3,6.3){$3$}
\put(6.3,6.3){$5$}
\put(7.3,6.3){$6$}
\put(8.3,6.3){$7$}
\put(9.3,6.3){$9$}

\put(0.3,5.3){$2$}
\put(1.3,5.3){$2$}
\put(2.3,5.3){$3$}
\put(3.3,5.3){$3$}
\put(4.3,5.3){$4$}
\put(5.3,5.3){$4$}
\put(6.3,5.3){$5$}

\put(0.3,4.3){$3$}
\put(1.3,4.3){$4$}
\put(2.3,4.3){$7$}
\put(3.3,4.3){$9$}

\put(0.3,3.3){$4$}
\put(1.3,3.3){$6$}
\put(2.3,3.3){$8$}

\put(0.3,2.3){$5$}
\put(1.3,2.3){$6$}

\put(0.3,1.3){$7$}

\put(0.3,0.3){$7$}

\end{picture}
\end{center}
\setlength{\unitlength}{1pt}

\noindent
The thicker line indicates the subtableau formed by
the entries in $\{1,2,3,4\}$.
\qed
\eex

Given such a supertableau $\La$,
for any $1\leqslant i\leqslant s\leqslant m$
denote by $\la_{si}$ the number of entries in row $i$
which do not exceed $s$. Furthermore, set $r=\la_{m1}$ and
for any $0\leqslant p\leqslant n$
and $1\leqslant j\leqslant r+p$
denote by $\la'_{r+p,j}$ the number of entries in column $j$
which do not exceed $m+p$.

The subtableau of the supertableau $\La$ occupied by the entries
in $\{1,\dots,m\}$ is a {\it column-strict\/} tableau
of the shape $\mu=(\la_{m1},\dots,\la_{mm})$.
The subtableau $\Tc$ of $\La$ occupied by the entries
in $\{m+1,\dots,m+n\}$ is a {\it row-strict\/} tableau
of the skew shape $\Ga_{\la}/\mu$.
We set
\ben
l_i=\la_i-i+1,\quad
l_{si}=\la_{si}-i+1,\quad
l^{\tss\prime}_{r+p,j}=\la'_{r+p,j}-j+1.
\een
For each $s=1,\dots,m+n$ we denote by $\om_s$ the number
of entries in $\La$ equal to $s$.

Our main theorem (Theorem~\ref{thm:baslla}) states that
each covariant representation $L(\la)$ admits a basis $\ze_{\La}$
parameterized by all supertableaux $\La$ of shape $\Ga_{\la}$.
Moreover, the action of the generators of the Lie superalgebra
$\gl_{m|n}$ in this basis is given by the formulas
\begin{align}\label{ekk}
E_{ss}\ts \ze^{}_{\La}&= \om_s\ts\ze^{}_{\La},
\\
\label{ekkp}
E_{s,s+1}\ts \ze^{}_{\Lambda}
&=\sum_{\La'} c_{\La\La'}\ts
\ze^{}_{\La'},\\
\label{ekkm}
E_{s+1,s}\ts \ze^{}_{\Lambda}
&=\sum_{\La'} d_{\La\La'}\ts
\ze^{}_{\La'},
\end{align}
where the sums in \eqref{ekkp} and \eqref{ekkm} are taken
over supertableaux $\La'$ obtained from $\La$
respectively by replacing
an entry $s+1$ by $s$ and by replacing
an entry $s$ by $s+1$. The coefficients $c_{\La\La'}$
and $d_{\La\La'}$ are found by explicit formulas
which depend on $s$, as well as on the row or column number of
$\La$ where the replacement occurs.
The symbols $\wedge_i$ in the denominators below indicate that
the $i$-th factor should be skipped.
For the values $s\ne m$ the coefficients
are given by
\ben
\bal
c_{\La\La'}&=-\frac{(l^{}_{si}-l^{}_{s+1,1})
\dots (l^{}_{si}-l^{}_{s+1,\tss s+1})}
{(l^{}_{si}-l^{}_{s1})\ldots\wedge_i\ldots(l^{}_{si}-l^{}_{ss})},\\
d_{\La\La'}&=\frac{(l^{}_{si}-l^{}_{s-1,1})
\dots (l^{}_{si}-l^{}_{s-1,\tss s-1})}
{(l^{}_{si}-l^{}_{s1})\ldots\wedge_i\ldots(l^{}_{si}-l^{}_{ss})},
\eal
\een
if $1\leqslant s\leqslant m-1$ and
the replacement occurs in row $i$, and by
\ben
\bal
c_{\La\La'}&=-\frac{(l^{\ts\prime}_{r+p,j}-l^{\ts\prime}_{r+p+1,1})
\dots (l^{\ts\prime}_{r+p,j}-l^{\ts\prime}_{r+p+1,\tss r+p+1})}
{(l^{\ts\prime}_{r+p,j}-l^{\ts\prime}_{r+p,1})\ldots
\wedge_j\ldots(l^{\ts\prime}_{r+p,j}-l^{\ts\prime}_{r+p,\tss r+p})},\\
d_{\La\La'}&=\frac{(l^{\ts\prime}_{r+p,j}-l^{\ts\prime}_{r+p-1,1})
\dots (l^{\ts\prime}_{r+p,j}-l^{\ts\prime}_{r+p-1,\tss r+p-1})}
{(l^{\ts\prime}_{r+p,j}-l^{\ts\prime}_{r+p,1})\ldots
\wedge_j\ldots(l^{\ts\prime}_{r+p,j}-l^{\ts\prime}_{r+p,\tss r+p})},
\eal
\een
if $s=m+p$ for $1\leqslant p\leqslant n-1$ and
the replacement occurs in column $j$.
In the case $s=m$ the coefficients $c_{\La\La'}$
and $d_{\La\La'}$ are given by more complicated
formulas; see Theorem~\ref{thm:baslla}. This
includes the case where the parameter
$r$ is changed to $r+1$ or $r-1$ which results in additional
factors in the expressions for $c_{\La\La'}$ and $d_{\La\La'}$.

These matrix element formulas can be interpreted
by using an equivalent combinatorial description
of the basis vectors. Namely, to each supertableau $\La$
of shape $\Ga_{\la}$ we can associate the pair of arrays
of nonnegative integers
$(\Uc,\Vc)$ of the form
\begin{align}
&\qquad\la^{}_{m1}\qquad\la^{}_{m2}
\qquad\qquad\cdots\qquad\qquad\la^{}_{mm}\non\\
&\qquad\qquad\la^{}_{m-1,1}\qquad\ \ \cdots\ \
\ \ \qquad\la^{}_{m-1,m-1}\non\\
\Uc=&\quad\qquad\qquad\cdots\qquad\cdots\qquad\cdots\non\\
&\quad\qquad\qquad\qquad\la^{}_{21}\qquad\la^{}_{22}\non\\
&\quad\qquad\qquad\qquad\qquad\la^{}_{11}  \non
\end{align}
and
\begin{align}
\la^{\prime}_{r+n,1}\quad\la^{\prime}_{r+n,2}\qquad\cdots
\qquad\qquad&\cdots\qquad\qquad\qquad\la^{\prime}_{r+n,r+n}\non\\
\Vc=\qquad\qquad\dddots\quad\dddots
\qquad\cdots\qquad\quad&\cdots\qquad\qquad\qquad\antiddots\non\\
\qquad\qquad\la^{\prime}_{r+1,1}\quad \la^{\prime}_{r+1,2}
\qquad\quad&\cdots\quad\qquad\la^{\prime}_{r+1,r+1}\non\\
\quad\qquad\qquad\qquad\la^{\prime}_{r,1}
\qquad\ \la^{\prime}_{r,2}
\qquad&\cdots\qquad\la^{\prime}_{r,r}\non
\end{align}

\medskip
\noindent
By the properties of the supertableau $\La$, both $\Uc$ and $\Vc$ are
{\it patterns\/} as the following {\it betweenness\/} (or
{\it interlacing\/})
conditions hold:
\ben
\la^{}_{k+1,i}\geqslant\la^{}_{ki}\geqslant\la^{}_{k+1,i+1}
\qquad\text{for}\quad
1\leqslant i\leqslant k\leqslant m-1
\een
and
\ben
\la^{\prime}_{r+p,\tss j}\geqslant\la^{\prime}_{r+p-1,\tss j}
\geqslant\la^{\prime}_{r+p,\tss j+1}
\qquad\text{for}\quad
p=1,\dots,n\qquad\text{and}\quad j=1,\dots,r+p-1.
\een

The basis elements $E_{ij}$ of $\gl_{m|n}$ with
$1\leqslant i,j\leqslant m$ span a subalgebra isomorphic to $\gl_m$.
The action of the elements of this subalgebra affects only the
pattern $\Uc$ leaving $\Vc$ unchanged, while the action
of the subalgebra isomorphic to $\gl_n$
which is spanned by the $E_{ij}$ with
$m+1\leqslant i,j\leqslant m+n$, affects only the
pattern $\Vc$ leaving $\Uc$ unchanged.
The above formulas for the action of the generators of both of these
subalgebras coincide with the classical Gelfand--Tsetlin formulas.
The action of the odd generators $E_{m,m+1}$ and $E_{m+1,m}$
affects both $\Uc$ and $\Vc$.

For any complex number $a$ the mapping
\ben
E_{ij}\mapsto E_{ij}+\de_{ij}(-1)^{\bi}a
\een
defines an automorphism of the universal enveloping algebra
$\U(\gl_{m|n})$. Twisting $L(\la)$ by such an automorphism
amounts to the shift $\la_i\mapsto \la_i+(-1)^{\bi}a$ of
the components of $\la$. Hence, the basis provided
by Theorem~\ref{thm:baslla} and the matrix element
formulas are also valid for any representation $L(\la)$
which is isomorphic to the composition of a covariant
representation with such an automorphism.

The tuple of complex numbers
$\om=(\om_1,\dots,\om_m\tss|\tss \om_{m+1},\dots,\om_{m+n})$
is a {\it weight\/} of $L(\la)$
if the subspace
\ben
L(\la)_{\om}=\{\eta\in L(\la)\ |\ E_{ii}\tss\eta=\om_i\tss\eta,
\quad i=1,\dots,m+n\}
\een
is nonzero.
Let $x_1,\dots,x_{m+n}$ be indeterminates. Then
the {\it character\/} of $L(\la)$ is the polynomial
\beql{chardef}
\ch L(\la)=\sum_{\om}\dim L(\la)_{\om}\ts
x^{\om_1}_1,\dots,x^{\om_{m+n}}_{m+n}.
\eeq
By Theorem~\ref{thm:baslla}, $\dim L(\la)_{\om}$ equals
the number of supertableaux $\La$ containing $\om_s$
entries equal to $s$ for each $s=1,\dots,m+n$.
Hence, we recover the formula for $\ch L(\la)$
originally obtained in \cite{br:hy} and \cite{s:ta};
see Corollary~\ref{cor:char} below.

\section{Mickelsson--Zhelobenko superalgebra
and Yangian}\label{sec:mz}
\setcounter{equation}{0}

We start by applying the standard methods of the
Mickelsson algebra theory developed by Zhelobenko
(see \cite{z:it, z:rr}) to the Lie superalgebra
$\gl_{m|n}$ and its natural subalgebra $\gl_m$.
The corresponding Mickelsson--Zhelobenko superalgebra
will be described in a way
similar to its even counterpart;
cf. \cite{m:yt}.

Let $\h$ denote the Cartan subalgebra of $\gl_m$
spanned by the elements $E_{11},\dots,E_{mm}$ and let
$\Rr(\h)$ denote
the field of fractions of the commutative algebra $\U(\h)$.
Consider the extension of the
universal enveloping algebra $\U(\gl_{m|n})$ defined by
\ben
\U'(\gl_{m|n})=\U(\gl_{m|n})\ot_{\U(\h)} \Rr(\h).
\een
Let $\Jr$
denote the left ideal of $\U'(\gl_{m|n})$ generated by
the elements $E_{ij}$ with $1\leqslant i<j\leqslant m$.
Then $\Jr$ is a two-sided ideal of the normalizer
\beql{normaz}
\Norm \Jr=\{x\in\U'(\gl_{m|n})\ |\ \Jr x\subseteq \Jr\}.
\eeq
The {\it Mickelsson--Zhelobenko superalgebra\/} $\Zr(\gl_{m|n},\gl_m)$
is defined as the quotient algebra
\ben
\Zr(\gl_{m|n},\gl_m)=\Norm \Jr/\Jr.
\een
Then $\Zr(\gl_{m|n},\gl_m)$ is a
superalgebra over $\CC$ and a
natural left and right $\Rr(\h)$-module.

Generators of the superalgebra $\Zr(\gl_{m|n},\gl_m)$
can be constructed by using the {\it extremal projector\/} $p$
for the Lie algebra $\gl_m$; see \cite{ast:po}.
The projector $p$ is an element of an algebra $\Fr(\gl_m)$
of formal series of elements of $\U(\gl_m)$
and can be defined as follows. The positive roots of
$\gl_m$ with respect to $\h$
are naturally enumerated by the pairs of indices $(i,j)$
such that $1\leqslant i<j\leqslant m$ so that the element
$E_{ij}$ is the corresponding root vector.
Call a linear ordering of the positive roots
{\it normal\/}
if any composite root lies between its components.
Set
\ben
p_{ij}=\sum_{k=0}^{\infty}(E_{ji})^k(E_{ij})^k\ts
\frac{(-1)^k}{k!\ts (h_i-h_j+1)\dots (h_i-h_j+k)},
\een
where $h_i=E_{ii}-i+1$. Then $p$ is given by the product
\ben
p=\prod_{i<j}\ts p_{ij}
\een
taken in any normal ordering on the pairs
$(i,j)$.
The formal series $p$ does not depend on the normal ordering
and has the properties
\ben
E_{ij}\ts p=p\tss E_{ji}=0\qquad\text{for}\quad
1\leqslant i<j\leqslant m.
\een
Moreover, $p$
satisfies the conditions $p^2=p$ and $p^*=p$, where
$x\mapsto x^*$ is the involutive anti-automorphism of the algebra
$\Fr(\gl_m)$ such that $(E_{ij})^*=E_{ji}$.

The extremal projector $p$ can be regarded as
a natural operator on the quotient space $\U'(\gl_{m|n})/\Jr$.
The Mickelsson--Zhelobenko superalgebra
$\Zr(\gl_{m|n},\gl_m)$
can be identified with the image of $\U'(\gl_{m|n})/\Jr$:
\ben
\Zr(\gl_{m|n},\gl_m)=p\big(\U'(\gl_{m|n})/\Jr\big).
\een
For $i=1,\dots,m$ and $a=m+1,\dots,m+n$ introduce
elements of $\Zr(\gl_{m|n},\gl_m)$ by
\begin{align}
z_{ia}&=pE_{ia}(h_i-h_1)\dots (h_i-h_{i-1}),
\non\\
z_{ai}&=pE_{ai}(h_i-h_{i+1})\dots (h_i-h_m).
\non
\end{align}
The $\ZZ_2$-grading on the superalgebra $\Zr(\gl_{m|n},\gl_m)$
is inherited from that of the superalgebra
$\U(\gl_{m|n})$ so that all
the elements $z_{ia}$, $z_{ai}$ are odd. Together with
the even elements $E_{ab}$ with $a,b\in\{m+1,\dots,m+n\}$
they generate the Mickelsson--Zhelobenko superalgebra
in the sense that monomials in the $z_{ia}$, $z_{ai}$ and $E_{ab}$
span $\Zr(\gl_{m|n},\gl_m)$ as a left (or right) $\Rr(\h)$-module.

The explicit formulas for the elements $z_{ia}$ and $z_{ai}$
(modulo $\Jr$)
have the form
\begin{align}\label{simae}
z_{ia}&=
\sum_{i>i_1>\dots>i_s\geqslant 1}
E_{ii_1}E_{i_1i_2}\dots E_{i_{s-1}i_s}E_{i_sa}
(h_i-h_{j_1})\dots (h_i-h_{j_r}),\\
\label{smaie}
z_{ai}&=
\sum_{i<i_1<\dots<i_s\leqslant m}
E_{i_1i}E_{i_2i_1}\dots E_{i_{s}i_{s-1}}E_{ai_s}
(h_i-h_{j_1})\dots (h_i-h_{j_r}),
\end{align}
where $s=0,1,\dots$ and $\{j_1,\dots,j_r\}$ is the complementary subset
to $\{i_1,\dots,i_s\}$ respectively in the set $\{1,\dots,i-1\}$
or $\{i+1,\dots,m\}$.

In the following proposition we let
the indices $i,j$ and $a,b,c$ run over
the sets $\{1,\dots,m\}$ and $\{m+1,\dots,m+n\}$, respectively.

\bpr\label{prop:relmz}
The following relations hold in $\Zr(\gl_{m|n},\gl_m)${\rm:}
\beql{Ess}
[E_{ab}, z_{ci}]=\delta_{bc}z_{ai},\qquad
[E_{ab}, z_{ic}]=-\delta_{ac}z_{ib}.
\eeq
Moreover, if $i\ne j$ then
\ben
z_{ai}z_{bj}=-z_{bj}z_{ai}\ts\frac{h_i-h_j+1}{h_i-h_j}+
z_{bi}z_{aj}\ts\frac{1}{h_i-h_j},
\een
while
\beql{siasja}
z_{ai}z_{bi}=-z_{bi}z_{ai}.
\eeq
Finally, if $i\ne j$ then
\beql{ziazbj}
z_{ia}z_{bj}=-z_{bj}z_{ia},
\eeq
while
\beq\label{siasbi}
z_{ia}z_{bi}=\big(\delta_{ba}(h_i+m-1)-E_{ba}\big)
\prod_{j=1,\ts j\ne i}^m(h_i-h_j-1)
-\sum_{j=1}^m z_{bj}z_{ja}
\prod_{k=1,\ts k\ne j}^m \frac{h_i-h_k-1}{h_j-h_k}.
\eeq
\epr

\bpf
All relations are verified by a standard calculation involving
the extremal projector $p$; cf. \cite{z:za}. They differ
only by signs from the relations in the Mickelsson--Zhelobenko
algebra $\Zr(\gl_{m+n},\gl_m)$; see \cite{m:yt}.
\epf

For all $a,b\in\{m+1,\dots,m+n\}$
introduce the polynomials
$Z_{ab}(u)$
in a variable $u$ with coefficients
in $\Zr(\gl_{m+n},\gl_m)$ by
\beql{bzab}
Z_{ab}(u)=\Big(\de_{ab}(u+m)+E_{ab}\Big)
\prod_{i=1}^{m}(u-h_i)
+\sum_{i=1}^{m}z_{ai}z_{ib}
\prod_{j=1,\ts j\ne  i}^{m}\frac{u-h_j}{h_i-h_j}.
\eeq
We will need to evaluate $Z_{ab}(u)$ at
$u=h$, where $h$ is an element of $\U(\h)$.
In order to make such evaluations
unambiguous, we will always assume that
the coefficients of the polynomial
are written to the left of the powers of $u$.
In particular, we have
\beql{zabhi}
Z_{ab}(h_i)=z_{ai}z_{ib},\qquad i\in\{1,\dots,m\}.
\eeq
Moreover, the relation \eqref{siasbi} implies that
\beql{zabhimo}
Z_{ab}(h_i-1)=-z_{ib}z_{ai},\qquad i\in\{1,\dots,m\}.
\eeq
This implies an alternative formula for the polynomials $Z_{ab}(u)$,
\beql{bzabinv}
Z_{ab}(u)=\Big(\de_{ab}u+E_{ab}\Big)
\prod_{i=1}^{m}(u-h_i+1)
-\sum_{i=1}^{m}z_{ib}z_{ai}
\prod_{j=1,\ts j\ne  i}^{m}\frac{u-h_j+1}{h_i-h_j}.
\eeq

The mapping $E_{ij}\mapsto (E_{ij})^*=E_{ji}$ defines
an involutive
anti-automorphism of the universal enveloping algebra
$\U(\gl_{m|n})$. We will use the same notation
for its natural extension to $\U'(\gl_{m|n})$
and to the Mickelsson--Zhelobenko superalgebra $\Zr(\gl_{m+n},\gl_m)$.
The elements of $\Rr(\h)$ are fixed points of
this anti-automorphism, and
the images of the generators of the superalgebra
are easy to calculate. It is also easy to verify directly that
the formulas in the next proposition define
an involutive anti-automorphism of $\Zr(\gl_{m+n},\gl_m)$.

\bpr\label{prop:imsig}
For any $i\in\{1,\dots,m\}$ and $a,b\in\{m+1,\dots,m+n\}$ we have
\ben
\bal
(z_{ia})^*&= z_{ai}\ts
\frac{(h_i-h_1-1)\dots(h_i-h_{i-1}-1)}{(h_i-h_{i+1})\dots(h_i-h_m)}\\
(z_{ai})^*&= z_{ia}\ts
\frac{(h_i-h_{i+1}+1)\dots(h_i-h_m+1)}{(h_i-h_1)\dots(h_i-h_{i-1})}
\eal
\een
and $(E_{ab})^*=E_{ba}$. Moreover, $\big(Z_{ab}(u)\big)^*=Z_{ba}(u)$.
\qed
\epr

\bpr\label{prop:relymz}
For any $i\in\{1,\dots,m\}$ and $a,b,c\in\{m+1,\dots,m+n\}$
the following relations hold in the Mickelsson--Zhelobenko
superalgebra $\Z(\gl_{m|n},\gl_m)$,
\begin{align}\label{zabuzci}
Z_{ab}(u)z_{ci}&=z_{ci}Z_{ab}(u)\ts\frac{u-h_i+1}{u-h_i}+
z_{ai}Z_{cb}(u)\ts\frac{1}{u-h_i}\\
\label{zabuzciinv}
Z_{ab}(u)z_{ci}&=z_{ci}Z_{ab}(u)\ts\frac{u-h_i+2}{u-h_i+1}+
Z_{cb}(u)z_{ai}\ts\frac{1}{u-h_i+1}.
\end{align}
\epr

\bpf
The first relation is verified by a straightforward calculation
with the use of the relations in $\Z(\gl_{m|n},\gl_m)$
given in Proposition~\ref{prop:relmz}. The second
relation is obtained by writing another
form of \eqref{zabuzci} with $a$ and $c$
swapped and solving the system of two equations simultaneously
for the unknowns $z_{ci}Z_{ab}(u)$ and $z_{ai}Z_{cb}(u)$.
\epf

Note that analogous relations involving the raising
operators $z_{ia}$ can be obtained by applying the anti-automorphism
of Proposition~\ref{prop:imsig} to the relations
of Proposition~\ref{prop:relymz}.

\medskip

Recall that the {\it Yangian for $\gl_n$\/} is
a unital associative algebra
$\Y(\gl_n)$ over $\CC$
with countably many generators $t_{ij}^{(1)},\ t_{ij}^{(2)},\dots$
where $i$ and $j$ run over the set $\{1,\dots,n\}$.
The defining relations of $\Y(\gl_n)$
have the form
\beql{defyang}
[t^{(r+1)}_{ij}, t^{(s)}_{kl}]-[t^{(r)}_{ij}, t^{(s+1)}_{kl}]=
t^{(r)}_{kj} t^{(s)}_{il}-t^{(s)}_{kj} t^{(r)}_{il},
\eeq
where $r,s\geqslant 0$ and $t^{(0)}_{ij}:=\delta_{ij}$.
Using the formal generating
series
\ben
t_{ij}(u) = \delta_{ij} + t^{(1)}_{ij} u^{-1} + t^{(2)}_{ij}u^{-2} +
\dots
\een
the defining relations
\eqref{defyang} can be written in the equivalent form
\ben
(u-v)\ts [t_{ij}(u),t_{kl}(v)]=
t_{kj}(u)\ts t_{il}(v)-t_{kj}(v)\ts t_{il}(u).
\een
A detailed exposition of the algebraic structure
and representation theory of the Yangian can be found in
\cite{m:yc}.

\bpr\label{prop:ymzahom}
The mapping
\ben
t_{ij}(u)\mapsto Z_{m+i,m+j}(u)\ts\frac{1}{(u+m)(u-h_1)\dots(u-h_m)}
\een
defines an algebra homomorphism
\beql{yamzhom}
\vp:\Y(\gl_n)\to \Zr(\gl_{m|n},\gl_m).
\eeq
\epr

\bpf
This is a super-analogue of the homomorphism from
the algebra $\Y(\gl_n)$ to the Mickelsson--Zhelobenko
algebra $\Zr(\gl_{m+n},\gl_m)$; see \cite[Theorem~3.1]{m:yt}.
Similar to the proof of that result, we note first
that the normalizer $\Norm\tss\Jr$ defined in \eqref{normaz}
contains the centralizer $\U(\gl_{m|n})^{\gl_m}$ of $\gl_m$
as a natural subalgebra. We will then
obtain the homomorphism \eqref{yamzhom} as the composition
of the super version of the {\it Olshanski homomorphism\/}
$\Y(\gl_n)\to\U(\gl_{m|n})^{\gl_m}$ and the natural homomorphism
$\U(\gl_{m|n})^{\gl_m}\to \Z(\gl_{m|n},\gl_m)$.

In more detail, let $E$ denote the $m\times m$ matrix whose
$(i,j)$ entry is the basis element $E_{ij}$ of $\gl_{m}$.
Then the mapping $\psi:\Y(\gl_n)\to\U(\gl_{m|n})$ given by
\beql{repdef}
\bal
t_{ij}^{(1)}&\mapsto E_{m+i,m+j},\\
t_{ij}^{(r)}&\mapsto \sum_{k,l=1}^m
E_{m+i,k}(E^{\tss r-2})_{kl}E_{l,m+j},\qquad r\geqslant 2,
\eal
\eeq
defines an algebra homomorphism. It can be verified
directly that the images of the generators $t_{ij}^{(r)}$ satisfy
the defining relations of the Yangian. Alternatively,
we can write this map as the composition of an
embedding of $\Y(\gl_n)$ into the Yangian $\Y(\gl_{m|n})$
and the evaluation homomorphism $\Y(\gl_{m|n})\to\U(\gl_{m|n})$;
see \cite[Formulas~(7) and (12)]{g:gd}.

Furthermore, it is easy to verify that each element
$\psi(t_{ij}^{(r)})$ commutes with the elements
$E_{kl}$, $1\leqslant k,l\leqslant m$,
so that
the image of the homomorphism $\psi$ is contained
in the centralizer $\U(\gl_{m|n})^{\gl_m}$.
The next step is to find the images of the elements
$\psi(t_{ij}^{(r)})$ in the Mickelsson--Zhelobenko
superalgebra by using calculations similar to the
even case; cf. \cite[Sec.~3]{m:yt} and \cite[Sec.~9.3]{m:yc}.
The final formula for the images of the coefficients of the series
$t_{ij}(u)$ is obtained by twisting the homomorphism $\psi$ by
the shift automorphism of the Yangian sending $t_{ij}(u)$
to $t_{ij}(u+m)$.
\epf

\section{Yangian action on the multiplicity space}\label{sec:ya}
\setcounter{equation}{0}

For any tuple of complex numbers
$\la=(\la_1,\dots,\la_m\ts|\ts\la_{m+1},\dots,\la_{m+n})$
satisfying the conditions
\eqref{findim} consider the corresponding
finite-dimensional irreducible representation
$L(\la)$ of $\gl_{m|n}$. Denote by
$L(\la)^+$ the subspace of
$\gl_m$-highest vectors in $L(\la)$:
\ben
L(\la)^+=\{\eta\in L(\la)\ts|\ts
E_{ij}\eta=0\ts
\qquad\text{for}\quad 1\leqslant i<j\leqslant m\}.
\een
Given an $m$-tuple of complex numbers $\mu=(\mu_1,\dots,\mu_m)$
such that $\mu_i-\mu_{i+1}\in\ZZ_+$ for all $i$
we denote by $L(\la)^+_{\mu}$ the corresponding weight subspace of
$L(\la)^+$:
\ben
L(\la)^+_{\mu}=\{\eta\in L(\la)^+\ts|\ts E_{ii}\eta=
\mu_i\ts\eta
\qquad\text{for}\quad i=1,\dots,m\}.
\een
We have the weight space decomposition
\ben
L(\la)^+=\bigoplus_{\mu} L(\la)^+_{\mu}.
\een

It follows from
the super-extension of the
general results of the Mickelsson algebra
theory (\cite{z:it}, \cite[Theorem~4.3.8]{z:rr}),
that given a total weight-consistent order on the set of
elements $E_{ba}$ with $m+1\leqslant a<b\leqslant m+n$
and $z_{ai}$ with $a=m+1,\dots,m+n$ and $i=1,\dots,m$,
the subspace $L(\la)^+$ is spanned by the
vectors $M\ze$, where $M$ runs over
the set of ordered monomials in these elements
and $\ze$ is the highest vector of $L(\la)$.
This implies that
the multiplicity space $L(\la)^+_{\mu}$ is nonzero
only if all components of the
$m$-tuple $\mu=(\mu_1,\dots,\mu_m)$
satisfy the inequalities
$0\leqslant \la_i-\mu_i\leqslant n$.
However, we will not need to rely on this result as its
independent proof will follow from the explicit construction
of the representation $L(\la)$. Namely,
assuming in addition that $\mu$ is a partition, we will construct
a basis of each space $L(\la)^+_{\mu}$ and use these bases to
produce
a $\gl_{m|n}$-submodule $K$ of $L(\la)$. Since $L(\la)$
is irreducible, we will conclude that $K=L(\la)$ so that
the space $L(\la)^+_{\mu}$ is nonzero if and only if
$\mu$ satisfies the above conditions.

The dimension of $L(\la)^+_{\mu}$ coincides with
the multiplicity of the $\gl_m$-module
$L'(\mu)$ in
the restriction of
$L(\la)$ to $\gl_m$.
The multiplicity
space $L(\la)^+_{\mu}$ is
a representation of the centralizer
$\U(\gl_{m|n})^{\gl_m}$.
Therefore, using the Olshanski
homomorphism $\Y(\gl_n)\to \U(\gl_{m|n})^{\gl_m}$
we can equip $L(\la)^+_{\mu}$ with an action of the Yangian.
As with the case of the skew representations
of the Yangian associated with the pair of Lie algebras
$\gl_{m+n}$ and $\gl_m$ (see e.g. \cite[Sec.~8.5]{m:yc})
one can show by extending
the arguments to the super case
that if the space $L(\la)^+_{\mu}$ is nonzero then
the resulting representation
of $\Y(\gl_n)$ in $L(\la)^+_{\mu}$ is irreducible.
We do not bring full details here as this would
lengthen the paper significantly and also because
the irreducibility of the $\Y(\gl_n)$-module
$L(\la)^+_{\mu}$ could be established in a more direct way
by using the character formula for $L(\la)$
obtained in \cite{br:hy} and \cite{s:ta}; see also
Corollary~\ref{cor:char} below.
As we want to demonstrate how this formula follows
from the basis construction, in our exposition
we will still rely on the properties of
the super-version of the Olshanski homomorphism.

\subsection{Covariant representations}\label{subsec:ya}

We will now suppose that all components $\la_i$ of $\la$
are nonnegative integers
and the number $\ell$
of nonzero components among $\la_{m+1},\dots,\la_{m+n}$
does not exceed $\la_m$.
The cyclic $\U(\gl_m)$-span of every
nonzero element of the multiplicity
space $L(\la)^+_{\mu}$ is
a finite-dimensional representation of $\gl_m$
isomorphic to $L'(\mu)$.

The elements $z_{ia}$ and $z_{ai}$ which are given by
explicit formulas \eqref{simae} and \eqref{smaie}
preserve the subspace of $\gl_m$-highest vectors in $L(\la)$.
Moreover, they raise and lower the $\gl_m$-weights,
respectively:
\ben
z_{ia}:L(\la)^+_{\mu}\to L(\la)^+_{\mu+\de_i},\qquad
z_{ai}:L(\la)^+_{\mu}\to L(\la)^+_{\mu-\de_i},
\een
where $\mu\pm\de_i$ is obtained from $\mu$ by replacing
the component $\mu_i$ by $\mu_i\pm 1$. We will call
the $z_{ia}$ and $z_{ai}$ the {\it raising\/} and
{\it lowering operators\/}, respectively.

Suppose that $\mu=(\mu_1,\dots,\mu_m)$ is a partition
satisfying the conditions $0\leqslant \la_i-\mu_i\leqslant n$ for
all $i=1,\dots,m$.
Introduce the element $\ze_{\mu}\in L(\la)^+_{\mu}$ by
\ben
\ze_{\mu}=\prod_{j=1}^m\big(z^{}_{m+\la_j-\mu_j,j}\dots
z^{}_{m+2,j}\tss z^{}_{m+1,j}\big)\ts\ze,
\een
with the product taken in the increasing order of $j$.

\bpr\label{prop:highvect}
Under
the action of $\Y(\gl_n)$ on $L(\la)^+_{\mu}$ we have
\ben
t_{ij}(u)\tss \ze_{\mu}=0\qquad\text{for}\quad
1\leqslant i<j\leqslant n
\een
and
\ben
t_{pp}(u)\tss\ze_{\mu}=\la_{m+p}(u)\tss\ze_{\mu}
\qquad\text{for}\quad p=1,\dots,n,
\een
where
\beql{lamk}
\la_{m+p}(u)=\frac{u+\la_{m+p}+m}{u+m}\tss
\prod_{i=1,\ \la_i-\mu_i\geqslant p}^m
\frac{u-\mu_i+i}{u-\mu_i+i-1}.
\eeq
\epr

\bpf
Using Proposition~\ref{prop:ymzahom}
we can rewrite the statements in terms of the action
of the coefficients of the polynomials $Z_{ab}(u)$.
Each element $h_i$ acts in $L(\la)^+_{\mu}$ as
multiplication by the scalar $\si_i=\mu_i-i+1$ so that
the relations can be written as
\beql{annihzemu}
Z_{ab}(u)\tss \ze_{\mu}=0\qquad\text{for}\quad
m+1\leqslant a<b \leqslant m+n
\eeq
and
\beql{eigenzemu}
Z_{m+p,m+p}(u)\tss\ze_{\mu}=
(u+\la_{m+p}+m)\tss
\prod_{j=1,\ \la_j-\mu_j\geqslant p}^m
(u-\si_j+1)\prod_{j=1,\ \la_j-\mu_j< p}^m(u-\si_j)\ts\ze_{\mu}
\eeq
for $p=1,\dots,n$.
We will prove them by induction on the sum
$\sum_{j=1}^m(\la_j-\mu_j)$. If the sum is zero, then
$\ze_{\mu}=\ze$ is the highest vector of $L(\la)$.
Then the relations follow from \eqref{bzab} as $z_{ia}\ts\ze=0$
for all $i$ and $a$. Now let $i$ be the minimum index such that
$\la_i-\mu_i>0$. Then $\ze_{\mu}=z_{ci}\tss \ze_{\mu+\de_i}$
with $c=m+\la_i-\mu_i$. Suppose that $a<b$. If
$a\geqslant c$, then we use \eqref{zabuzci}
to write
\beql{zabzvec}
Z_{ab}(u)\tss z_{ci}\ts \ze_{\mu+\de_i}
=z_{ci}Z_{ab}(u)\ts\frac{u-h_i+1}{u-h_i}\tss \ze_{\mu+\de_i}+
z_{ai}Z_{cb}(u)\ts\frac{1}{u-h_i}\tss \ze_{\mu+\de_i}.
\eeq
If $a<c$, then using \eqref{zabuzciinv} we get
\beql{zabzvecin}
Z_{ab}(u)\ts z_{ci}\ts \ze_{\mu+\de_i}
=z_{ci}Z_{ab}(u)\ts\frac{u-h_i+2}{u-h_i+1}\ts \ze_{\mu+\de_i}+
Z_{cb}(u)z_{ai}\ts\frac{1}{u-h_i+1}\ts \ze_{\mu+\de_i}.
\eeq
Note that $z_{ai}^2=0$ by \eqref{siasja}, so that
$z_{ai}\tss \ze_{\mu+\de_i}=0$ if $a<c$.
Therefore,
applying the induction hypothesis
we derive \eqref{annihzemu}.

Now let $a=b=m+p$. Applying \eqref{zabzvec} and \eqref{zabzvecin}
together with \eqref{annihzemu}, we get
\begin{alignat}{2}
Z_{aa}(u)\tss z_{ci}\ts \ze_{\mu+\de_i}
&=z_{ci}Z_{aa}(u)\ts\frac{u-\si_i}{u-\si_i-1}\ts \ze_{\mu+\de_i},
\qquad&&\text{if}\quad\la_i-\mu_i<p,
\non\\
Z_{aa}(u)\tss z_{ci}\ts \ze_{\mu+\de_i}
&=z_{ci}Z_{aa}(u)\ts\frac{u-\si_i+1}{u-\si_i}\ts \ze_{\mu+\de_i}
\qquad&&\text{if}\quad\la_i-\mu_i>p,
\non\\
Z_{aa}(u)\tss z_{ai}\ts \ze_{\mu+\de_i}
&=z_{ai}Z_{aa}(u)\ts\frac{u-\si_i+1}{u-\si_i-1}\ts \ze_{\mu+\de_i}
\qquad&&\text{if}\quad\la_i-\mu_i=p.
\non
\end{alignat}
This proves \eqref{eigenzemu}.
\epf

Recall the notation $l_j=\la_j-j+1$
for $j=1,\dots,m$. Fix an index
$i\in\{1,\dots,m\}$ and set $k=\la_i-\mu_i$.

\bco\label{cor:zact}
If $k=0$ then
$z_{i,m+p}\ts\ze_{\tss\mu}=0$ for all $p=1,\dots,n$.
If $k\geqslant 1$ then
\begin{align}
z_{i,m+k}\ts\ze_{\tss\mu}&=(\si_i+\la_{m+k}+m)
\ts\prod_{j=1}^{i-1}
(-1)^{\la_j-\mu_j}\ts(\si_i-l_j)
\non\\
\label{zimze}
{}&\times{}\prod_{j=i+1,\ \la_j-\mu_j\geqslant k}^{m}
(\si_i-\si_j+1)\ts
\prod_{j=i+1,\ \la_j-\mu_j< k}^{m}
(\si_i-\si_j)\ts\ze_{\tss\mu+\de_i}.
\end{align}
Moreover, for $0\leqslant k\leqslant n-1$ we have
\beql{zmze}
z_{m+k+1,i}\ts\ze_{\tss\mu}=
\prod_{j=1}^{i-1}
\frac{(-1)^{\la_j-\mu_j}}
{\si_i-l_j-1}\ts
\prod_{j=1,\ \la_j-\mu_j>k}^{i-1}
(\si_i-\si_j)\ts\prod_{j=1,\ \la_j-\mu_j\leqslant k}^{i-1}
(\si_i-\si_j-1)\ts
\ze_{\tss\mu-\de_i}.
\eeq
\eco

\bpf
Since $z_{i,m+p}\ts\ze=0$,
the first claim is immediate from \eqref{ziazbj}. Furthermore,
if $k\geqslant 1$, then
applying \eqref{ziazbj}, we get
\ben
z_{i,m+k}\ts\ze_{\tss\mu}
=\prod_{j=1}^{i-1}(-1)^{\la_j-\mu_j}
\prod_{j=1}^{i-1}\big(z^{}_{m+\la_j-\mu_j,j}\dots
z^{}_{m+2,j}\tss z^{}_{m+1,j}\big)\ts
z_{i,m+k}\ts\ze^{}_{\wt\mu},
\een
where $\wt\mu=(\la_1,\dots,\la_{i-1},\mu_i,\dots,\mu_m)$.
By \eqref{zabhimo}, we have
$z_{i,m+k}\ts z_{m+k,i}=-Z_{m+k,m+k}(h_i-1)$ so that
\ben
z_{i,m+k}\ts\ze^{}_{\wt\mu}
=z_{i,m+k}\ts z_{m+k,i}\ts \ze^{}_{\wt\mu+\de_i}
=-Z_{m+k,m+k}(\si_i)\ts\ze^{}_{\wt\mu+\de_i},
\een
and \eqref{zimze} follows from \eqref{eigenzemu}.

Now apply $z_{m+k,i}$ to both sides of \eqref{zimze}.
Using \eqref{zabhi} and \eqref{eigenzemu},
for the left hand side we obtain
\ben
\bal
z_{m+k,i}\ts z_{i,m+k}\ts\ze_{\mu}&=
Z_{m+k,m+k}(h_i)\ts\ze_{\mu}=Z_{m+k,m+k}(\si_i)\ts\ze_{\mu}\\
{}&=(\si_i+\la_{m+k}+m)\tss
\prod_{j=1,\ \la_j-\mu_j\geqslant k}^m
(\si_i-\si_j+1)\prod_{j=1,\ \la_j-\mu_j< k}^m
(\si_i-\si_j)\ts\ze_{\mu}.
\eal
\een
Rewriting the resulting relation for $\mu-\de_i$ instead of $\mu$,
we get \eqref{zmze}.
\epf

\bco\label{cor:nonzero}
Suppose that $\mu=(\mu_1,\dots,\mu_m)$ is a partition
such that the conditions $0\leqslant \la_i-\mu_i\leqslant n$ hold for
all $i=1,\dots,m$.
Then the vector $\ze_{\mu}$ is nonzero.
\eco

\bpf
As in the proof of Proposition~\ref{prop:highvect},
we argue by induction on the sum
$\sum_{j=1}^m(\la_j-\mu_j)$ and suppose that
$\ze_{\mu}=z_{m+k,i}\tss \ze_{\mu+\de_i}$
with $k=\la_i-\mu_i$, where $i$ is the minimum index such that
$\la_i-\mu_i>0$. Then $\si_j=l_j$ for $j=1,\dots,i-1$ and
using \eqref{zimze} we obtain
\ben
z_{i,m+k}\ts\ze_{\mu}=(\si_i+\la_{m+k}+m)\tss
\prod_{j=1,\ \la_j-\mu_j\geqslant k}^m
(\si_i-\si_j+1)\prod_{j=1,\ \la_j-\mu_j<k}^m
(\si_i-\si_j)\ts\ze_{\mu+\de_i}.
\een
We have
$\si_i+\la_{m+k}+m\geqslant\si_i+m>0$. Moreover,
$\si_1>\dots>\si_m$ so that the only factor in the coefficient
of $\ze_{\mu+\de_i}$ which could be equal to zero is
$\si_i-\si_{i-1}+1$. However, in this case
$\mu_i=\mu_{i-1}=\la_{i-1}$ which contradicts the assumption
$\la_i-\mu_i>0$. Thus, the coefficient is nonzero. By
the induction hypothesis, the vector $\ze_{\mu+\de_i}$ is nonzero
and hence so is $\ze_{\mu}$.
\epf

We will keep
the assumptions of Proposition~\ref{prop:highvect}
and Corollary~\ref{cor:nonzero}.
Due to these statements,
the irreducible $\Y(\gl_n)$-module
$L(\la)^+_{\mu}$ is a highest
weight representation whose
highest weight is the $n$-tuple
$(\la_{m+1}(u),\dots,\la_{m+n}(u))$. The corresponding {\it Drinfeld
polynomials\/} $P_1(u),\dots,P_{n-1}(u)$ are monic
polynomials in $u$ which are defined by
the relations
\beql{hwdp}
\frac{\la_{m+k}(u)}{\la_{m+k+1}(u)}=\frac{P_k(u+1)}{P_k(u)}
\eeq
for $k=1,\dots,n-1$; see e.g. \cite[Sec.~3.2]{m:yc}.
In Sec.~\ref{sec:mt} we associated a Young diagram
$\Ga_{\la}$ to each covariant highest weight $\la$.
To give formulas for the polynomials $P_k(u)$, consider
the skew diagram $\Ga_{\la}/\mu$ obtained from $\Ga_{\la}$
by removing the first $\mu_i$ boxes in row $i$ for each
$i=1,\dots,m$. By the {\it content\/} of any box
$\al$ of $\Ga_{\la}/\mu$ we will mean the number $c(\al)=j-i$
if $\al$ is the intersection of row $i$ and
column $j$ of the diagram.

\bco\label{cor:dripol}
The Drinfeld polynomials associated with the
$\Y(\gl_n)$-module
$L(\la)^+_{\mu}$ are given by
the formulas
\ben
P_k(u)=\prod_{\al} \big(u-c(\al)\big),
\een
where $\al$ runs over the leftmost boxes
of the rows of length $k$ in
the diagram $\Ga_{\la}/\mu$.
\eco

\bpf
The formulas of Proposition~\ref{prop:highvect} imply that
\ben
P_k(u)=\prod_{i=1,\ \la_i-\mu_i= k}^m(u-\si_i)
\prod_{j=1}^{\la_{m+k}-\la_{m+k+1}}(u+\la_{m+k+1}+m+j-1).
\een
Writing the factors in terms of the contents of
the boxes of the diagram $\Ga_{\la}/\mu$ gives
the desired formulas.
\epf

Now we introduce some parameters of the diagram conjugate
to $\Ga_{\la}/\mu$. Set $r=\mu_1$ and let
$\mu'=(\mu'_1,\dots,\mu'_r)$ be the diagram conjugate to
$\mu$ so that $\mu'_j$ equals the number of boxes in column
$j$ of $\mu$. Furthermore, set $\la'=(\la'_1,\dots,\la'_{r+n})$,
where $\la'_j$ equals the number of boxes in column
$j$ of the diagram $\Ga_{\la}$.

Consider the general linear Lie
algebra $\gl_{r+n}$ and its natural subalgebra $\gl_r$.
Denote by $\BL(\la')$ the finite-dimensional irreducible
representation of $\gl_{r+n}$ with the highest weight $\la'$.
The subspace $\BL(\la')^+_{\mu'}$ of $\gl_r$-highest vectors
in $\BL(\la')$ of weight $\mu'$ is equipped with
a structure of irreducible representation of
the Yangian $\Y(\gl_n)$. This is a
{\it skew representation\/} of $\Y(\gl_n)$; see e.g.
\cite[Sec. 8.5]{m:yc} for a more detailed description
of these representations. To define the action,
consider the homomorphism $\Y(\gl_n)\to\U(\gl_{r+n})$
defined by
\beql{actskew}
\bal
t_{ij}^{(1)}&\mapsto \BE_{r+i,r+j},\\
t_{ij}^{(p)}&\mapsto (-1)^{p-1}\sum_{k,l=1}^r
\BE_{r+i,k}({\wh E}^{\ts p-2})_{kl}\BE_{l,r+j},\qquad p\geqslant 2,
\eal
\eeq
where ${\wh E}$ denotes the $r\times r$ matrix whose
$(i,j)$ entry is the basis element $\BE_{ij}$ of $\gl_{r}$.
The image of the homomorphism is contained
in the centralizer $\U(\gl_{r+n})^{\gl_r}$ which allows to define
the Yangian action on the vector space $\BL(\la')^+_{\mu'}$
via this homomorphism. We will work with
the twisted action of the Yangian on this space which is obtained
by taking its composition with the automorphism
sending $t_{ij}(u)$
to $t_{ij}(u-r)$.

\bth\label{thm:ident}
The $\Y(\gl_n)$-modules $L(\la)^+_{\mu}$ and
$\BL(\la')^+_{\mu'}$ are isomorphic.
\eth

\bpf
The Drinfeld polynomials $\BP_k(u)$ of the skew
$\Y(\gl_n)$-module
$\BL(\la')^+_{\mu'}$ were first calculated in \cite{nt:ry}
for a slightly different action of the Yangian.
In our setting they can be written in the form
\ben
\BP_k(u)=\prod_{\al} \big(u+c(\al)\big),\qquad k=1,\dots,n-1,
\een
where $\al$ runs over the top boxes
of the columns of height $k$ in
the diagram $\la'/\mu'$; see also
\cite[Sec.~8.5]{m:yc}. Hence, by Corollary~\ref{cor:dripol}
the Drinfeld polynomials
of $\BL(\la')^+_{\mu'}$ coincide with those of
the module $L(\la)^+_{\mu}$. We want to verify that
the highest weight of this module
is the same as the highest weight $\nu(u)=(\nu_1(u),\dots,\nu_n(u))$
of the $\Y(\gl_n)$-module $\BL(\la')^+_{\mu'}$.
We will use the formula for $\nu(u)$
given in \cite[Theorem~8.5.4]{m:yc}: the components $\nu_k(u)$
are found by
\ben
\nu_k(u)=\frac{(u+\nu_k^{(1)})\tss(u+\nu_k^{(2)}-1)\dots
(u+\nu_k^{(r+1)}-r)}
{(u+\mu'_1)\tss(u+\mu'_2-1)\dots
(u+\mu'_r-r+1)\tss (u-r)},
\een
where
\ben
\nu^{(s)}_k=\middd\{\mu'_{s-1},\mu'_s,\la'_{k+s-1}\},
\een
assuming $\mu'_0$ is sufficiently large,
$\mu'_{r+1}=0$, and $\middd\{a,b,c\}$
denotes the middle of the three integers. Since
the components of the highest weight and the Drinfeld polynomials
are related by \eqref{hwdp}, it will be sufficient to
demonstrate that $\nu_1(u)=\la_{m+1}(u)$.
Write $\mu=(\mu_1,\dots,\mu_m)$ in the form
\ben
\mu=(\underbrace{r,\dots,r}_{\mu'_r},\ts
\underbrace{r-1,\dots,r-1}_{\mu'_{r-1}-\mu'_r},\dots,
\underbrace{1,\dots,1}_{\mu'_{1}-\mu'_2},\ts
\underbrace{0,\dots,0}_{m-\mu'_1}).
\een
Now calculate $\la_{m+1}(u)$ by \eqref{lamk}.
The part of the product
\ben
\prod_{i=1,\ \la_i-\mu_i\geqslant 1}^m
\frac{u-\mu_i+i}{u-\mu_i+i-1}
\een
corresponding to the
values $\mu_i=s$ simplifies to the expression
\ben
\frac{u+\nu^{(s+1)}_1-s}{u+\mu'_{s+1}-s}
\een
for each $s\in\{1,\dots,r\}$. The same expression
equals the part of this product with $\mu_i=s=0$
multiplied by the first factor in \eqref{lamk}.
\epf

\subsection{Construction of the basis vectors}\label{subsec:cbv}

As we pointed out in the proof of Theorem~\ref{thm:ident},
the representation $\BL(\la')^+_{\mu'}$ of $\Y(\gl_n)$
admits a basis parameterized by column-strict tableaux
of shape $\la'/\mu'$. The explicit action of the Drinfeld
generators of the Yangian in this basis
was given in \cite{nt:ry} by using a combinatorially equivalent
description of the basis vectors in terms of
trapezium-like patterns. Using Theorem~\ref{thm:ident},
we get a basis of the $\Y(\gl_n)$-module $L(\la)^+_{\mu}$
parameterized by the row-strict tableaux of shape $\Ga_{\la}/\mu$
with entries in $\{m+1,\dots,m+n\}$, together with
an explicit action of $\Y(\gl_n)$ in this basis.
Furthermore, taking the Gelfand--Tsetlin basis of
each $\gl_m$-module $L'(\mu)$ we get a basis of
the direct sum $K$ of the vector spaces in \eqref{decomptp}
taken over the highest weights $\mu$ satisfying
the assumptions of Corollary~\ref{cor:nonzero}.
Our goal now is to obtain explicit formulas for the action
of generators of $\gl_{m|n}$ in this basis of $K$. This
will show that $K$ is stable under the action
and hence, due to
the irreducibility of $L(\la)$, that $K=L(\la)$.
The first step will be to present the basis vectors
explicitly in terms of the lowering operators.

Given elements
$a_1,\dots,\tss a_k$ and $b_1,\dots,\tss b_k$
of the set $\{m+1,\dots,m+n\}$, introduce
the corresponding {\it quantum minors\/}
\beql{quamin}
Z^{\ts a_1,\dots,\tss a_k}_{\ts b_1,\dots,\tss b_k}(u)
=\sum_{p\in \Sym_k} \sgn p\cdot
Z_{a_{1}b_{p(1)}}(u-k+1)\dots
Z_{a_{k}b_{p(k)}}(u),
\eeq
where the polynomials $Z_{ab}(u)$ are defined by \eqref{bzab}.
By Proposition~\ref{prop:ymzahom}, the polynomials
\eqref{quamin} inherit the properties of the respective quantum
minors $t^{\ts a_1,\dots,\tss a_k}_{\ts b_1,\dots,\tss b_k}(u)$
in the Yangian $\Y(\gl_n)$; see e.g. \cite[Sec.~1.6]{m:yc}.
In particular, the polynomials are skew-symmetric with respect
to permutations of the upper indices $a_i$ or
lower indices $b_i$, and they can be expressed in a form
alternative to \eqref{quamin}:
\beql{quaminhv}
Z^{\ts a_1,\dots,\tss a_k}_{\ts b_1,\dots,\tss b_k}(u)
=\sum_{p\in \Sym_k} \sgn p\cdot
Z_{a_{p(1)}b_{1}}(u)\dots
Z_{a_{p(k)}b_{k}}(u-k+1).
\eeq
Moreover, their images with respect
to the anti-automorphism of Proposition~\ref{prop:imsig}
are given by
\beql{qmininant}
\big(Z^{\ts a_1,\dots,\tss a_k}_{\ts b_1,\dots,\tss b_k}(u)\big)^*
=Z^{\ts b_1,\dots,\tss b_k}_{\ts a_1,\dots,\tss a_k}(u).
\eeq

\ble\label{lem:relqm}
The following relations hold in $\Z(\gl_{m|n},\gl_m)${\rm :}
for any $c\in\{a_1,\dots,a_k\}$ we have
\begin{align}\label{aux}
z_{ci}\tss Z^{\ts a_1,\dots,\tss a_k}_{\ts b_1,\dots,\tss b_k}(u)
&=Z^{\ts a_1,\dots,\tss a_k}_{\ts b_1,\dots,\tss b_k}(u)
\tss z_{ci}\ts\frac{u-h_i-k+1}{u-h_i+2},\\
\label{auxtwo}
z_{ic}\tss Z^{\ts b_1,\dots,\tss b_k}_{\ts a_1,\dots,\tss a_k}(u)
&=Z^{\ts b_1,\dots,\tss b_k}_{\ts a_1,\dots,\tss a_k}(u)
\tss z_{ic}\ts\frac{u-h_i+1}{u-h_i-k}.
\end{align}
\ele

\bpf
Due to Proposition~\ref{prop:imsig} and \eqref{qmininant},
the relations are equivalent. Therefore it
suffices to prove \eqref{aux}.
By the skew-symmetry of the quantum minors, we may assume
that $c=a_k$. We will argue by
induction on $k$. If $k=1$, then \eqref{aux} holds by
\eqref{zabuzci} with $a=c$. For $k\geqslant 2$ write
\ben
Z^{\ts a_1,\dots,\tss a_k}_{\ts b_1,\dots,\tss b_k}(u)
=\sum_{j=1}^k (-1)^{j-1}Z_{a_1b_j}(u-k+1)
Z^{\ts a_2,\dots,\tss a_k}_{\ts b_1,\dots,\tss
\hat b_j,\dots,\tss b_k}(u).
\een
Hence, using again \eqref{zabuzci} we get
\ben
\bal
z_{ci}\tss Z^{\ts a_1,\dots,\tss a_k}_{\ts b_1,\dots,\tss b_k}(u)
=\sum_{j=1}^k (-1)^{j-1}
\Big(&Z_{a^{}_1b_j}(u-k+1)\tss z_{ci}\tss\frac{u-h_i-k+1}{u-h_i-k+2}\\
{}-z_{a^{}_1i}\tss &Z_{cb_j}(u-k+1)\tss\frac{1}{u-h_i-k+2}\Big)\ts
Z^{\ts a_2,\dots,\tss a_k}_{\ts b_1,\dots,\tss
\hat b_j,\dots,\tss b_k}(u).
\eal
\een
Since
$Z^{\ts c,\tss a_2,\dots,\tss a_k}_{\ts b_1,\dots,\tss b_k}(u)=0$,
applying the induction hypothesis we arrive at \eqref{aux}.
\epf

For each $p=1,\dots,n$ set
\beql{defaku}
A_p(u)=Z^{\ts m+1,\dots,\tss m+p}_{\ts m+1,\dots,\tss m+p}(u)
\ts
Z^{\ts m+1,\dots,\tss m+p-1}_{\ts m+1,\dots,\tss m+p-1}(u)^{-1}
\prod_{i=2}^{m+1}\frac{1}{u-h_i-p+1},
\eeq
and for $p=1,\dots,n-1$ set
\ben
B_p(u)=-Z^{\ts m+1,\dots,\tss m+p}_{\ts m+1,\dots,
\tss m+p-1,\tss m+p+1}(u)
\ts
Z^{\ts m+1,\dots,\tss m+p+1}_{\ts m+1,\dots,\tss m+p+1}(u)^{-1}
\prod_{i=2}^{m+1}(u-h_i-p)
\een
and
\ben
C_p(u)=Z^{\ts m+1,\dots,\tss m+p-1,
\tss m+p+1}_{\ts m+1,\dots,\tss m+p}(u)
\ts
Z^{\ts m+1,\dots,\tss m+p-1}_{\ts m+1,\dots,\tss m+p-1}(u)^{-1}
\prod_{i=2}^{m+1}\frac{1}{u-h_i-p+1},
\een
where $h_{m+1}=-m$ and the second quantum minor in the formulas
for $A_1(u)$ and $C_1(u)$ is understood as being equal to $1$.

As before, we also regard the $Z_{ab}(u)$ as polynomials
in $u$ whose coefficients are operators
in $L(\la)^+_{\mu}$. We will see below
that the basis vectors of $L(\la)^+_{\mu}$ are eigenvectors
for the coefficients of all
quantum minors
$Z^{\ts m+1,\dots,\tss m+l}_{\ts m+1,\dots,\tss m+l}(u)$.
Therefore, the application of $A_p(u)$, $B_p(u)$ or $C_p(u)$
to a basis vector produces a linear combination
of the basis vectors whose coefficients are rational functions
in $u$.

Recalling the parametrization of the basis vectors of the skew
representations of the Yangian (see \cite[Sec.~8.5]{m:yc})
and using the isomorphism of Theorem~\ref{thm:ident}, we find
that the highest vector
$\ze_{\mu}$ of the $\Y(\gl_n)$-module corresponds
to the {\it initial\/} $\Ga_{\la}/\mu$-{\it tableau\/}
$\Tc^{\tss 0}$ which
is obtained by filling in the boxes of each row by
the consecutive numbers $m+1,m+2,\dots$ from left to right.
The entries $\la_{r+p,j}^0$ of the corresponding
pattern $\Vc^{\tss 0}$ are given by
\beql{lazero}
\la_{r+p,j}^0=\min\{\la'_j,\mu'_{j-p}\},
\eeq
where we assume that $\mu'_{i}$ is sufficiently large
for $i\leqslant 0$. Equivalently,
the parameters $\la^0_{r+p,j}$ of $\Tc^{\tss0}$
can be defined by first extending
$\Tc^{\tss0}$ to a supertableau
of shape $\Ga_{\la}$ by writing the entry $i$ in each box
of row $i$ of $\mu$ for $i=1,\dots,m$,
so that
$\la^0_{r+p,j}$ is the number of entries in column $j$
of this supertableau which do not exceed $m+p$.
We will use the notation
\beql{lzero}
l^{\tss0}_{r+p,j}=\la^0_{r+p,j}-j+1.
\eeq

Given an arbitrary row-strict $\Ga_{\la}/\mu$-tableau $\Tc$,
set
\beql{basvdef}
\ze^{}_{\Tc}=\prod_{(p,j)}
\Big(C_p(-l'_{r+p,j}-1)\dots C_p(-l^{\tss 0}_{r+p,j}+1)\tss
C_p(-l^{\tss 0}_{r+p,j})\Big)\ts\ze_{\mu},
\eeq
where the product is taken over the pairs $(p,j)$ with
$p=1,\dots,n-1$ and $j=1,\dots,r+p$ in the order
\begin{multline}
(n-1,1),\dots,(1,1),\ (n-1,2),\dots,(1,2),\ \dots,
\ (n-1,r+1),\dots,(1,r+1),\\
(n-1,r+2),\dots,(2,r+2),\ (n-1,r+3),\dots,(3,r+3),\ \dots,\\
(n-1,r+n-2),(n-2,r+n-2),\ (n-1,r+n-1).
\non
\end{multline}

\bpr\label{prop:basskew}
All evaluations of $C_p(u)$ involved in the expression
\eqref{basvdef} are well-defined. The vectors
$\ze^{}_{\Tc}$ parameterized by the row-strict tableaux $\Tc$
form a basis of $L(\la)^+_{\mu}$. Moreover,
the action of the generators of the Lie subalgebra $\gl_n$
in this basis is given
by the formulas \eqref{ekk}, \eqref{ekkp} and
\eqref{ekkm} with $s\geqslant m+1$.
\epr

\bpf
The rational function $C_p(u)$ coincides with the image
of the series
\beql{sert}
t^{\ts 1,\dots,\tss p-1,\tss p+1}_{\ts 1,\dots,\tss p}(u)
\ts
t^{\ts 1,\dots,\tss p-1}_{\ts 1,\dots,\tss p-1}(u)^{-1}
\tss(u-h_1-p+1)
\eeq
under the homomorphism of Proposition~\ref{prop:ymzahom}.
Note that $h_1$ acts in $L(\la)^+_{\mu}$ as multiplication
by the scalar $\mu_1=r$. Now we find the image
of the series \eqref{sert} (with $h_1$ replaced
by $r$) in the skew representation
$\BL(\la')^+_{\mu'}$ of $\Y(\gl_n)$. The formulas
defining this representation can be written
in an equivalent form with the use of quantum minors
of the matrix $1+\BE u^{-1}$, where $\BE$ denotes the
$(r+n)\times(r+n)$ matrix whose $(i,j)$ entry is $\BE_{ij}$.
Namely, the representation is defined by
\ben
t_{ij}(u)\mapsto
\Big[\big(1+\BE\ts u^{-1}\big)^{1,\dots,r}_{1,\dots,r}\Big]^{-1}
\cdot\big(1+\BE\ts u^{-1}\big)^{1,\dots,
\tss r,\tss r+i}_{1,\dots,\tss r,\tss r+j}
\een
for $i,j\in\{1,\dots,n\}$.
Moreover, the images of the quantum minors
occurring in \eqref{sert} are then found by
\ben
\bal
t^{\ts 1,\dots,\tss p-1,\tss p+1}_{\ts 1,\dots,\tss p}(u)
&\mapsto
\Big[\big(1+\BE\ts u^{-1}\big)^{1,\dots,r}_{1,\dots,r}\Big]^{-1}
\cdot
\big(1+\BE\ts u^{-1}\big)^{1,\dots,\tss r,\tss r+1,
\dots,\tss r+p-1,\tss r+p+1}_{1,\dots,\tss r,\tss r+1,
\dots,\tss r+p},\\
t^{\ts 1,\dots,\tss p-1}_{\ts 1,\dots,\tss p-1}(u)
&\mapsto
\Big[\big(1+\BE\ts u^{-1}\big)^{1,\dots,r}_{1,\dots,r}\Big]^{-1}
\cdot
\big(1+\BE\ts u^{-1}\big)^{1,\dots,\tss r,\tss r+1,
\dots,\tss r+p-1}_{1,\dots,\tss r,\tss r+1,
\dots,\tss r+p-1};
\eal
\een
see e.g. \cite[Sec.~8.5]{m:yc} for proofs of these statements.
Hence, calculating the image of the series \eqref{sert}
in the skew representation
$\BL(\la')^+_{\mu'}$ we conclude that the rational
function $C_p(u)$, regarded as an operator
in $\BL(\la')^+_{\mu'}$ can be written as
\ben
C_p(u)=
\big(u+\BE\big)^{1,\dots,\tss r+p-1, \tss r+p+1}_{1,\dots,\tss r+p}
\cdot
\Big[\big(u+\BE\big)^{1,\dots,r+p-1}_{1,\dots,r+p-1}\Big]^{-1}.
\een
However, as was observed in \cite{nt:yg},
for an appropriate value of $u$ such an operator
takes a vector of the Gelfand--Tsetlin basis of $\BL(\la')^+_{\mu'}$
to another vector of this basis; see also \cite[Sec.~5.4]{m:yc}.
More precisely,
if $\Vc$ and $\Vc^{\tss-}$ are trapezium patterns of the form
described in Sec.~\ref{sec:mt} such that
$\Vc^{\tss-}$ is obtained from $\Vc$
by replacing an entry $\la'_{r+p,j}$
by $\la'_{r+p,j}-1$, then for the corresponding
basis vectors $\ze^{}_{\Vc}$ and $\ze^{}_{\Vc^{-}}$
of $\BL(\la')^+_{\mu'}$ we have
\ben
C_p(-l^{\tss\prime}_{r+p,j})\tss \ze^{}_{\Vc}=\ze^{}_{\Vc^{-}}.
\een
Applying the isomorphism of Theorem~\ref{thm:ident}
we conclude that
\ben
C_p(-l^{\tss\prime}_{r+p,j})\tss \ze^{}_{\Tc}=\ze^{}_{\Tc^{-}}
\een
in the representation $L(\la)^+_{\mu}$,
there $\Tc$ and $\Tc^{-}$ are the
row-strict tableaux corresponding
to the patterns $\Vc$ and $\Vc^{-}$, respectively, so that
$\Tc^{-}$ is obtained from
$\Tc$ by replacing an entry $m+p$ by $m+p+1$ in column $j$.
This shows that the vectors given by \eqref{basvdef}
are well-defined and they
form a basis of $L(\la)^+_{\mu}$.

Comparing the actions of Yangian in $L(\la)^+_{\mu}$ and
$\BL(\la')^+_{\mu'}$ given by
\eqref{repdef} and \eqref{actskew}, we can conclude that
the elements of the subalgebra $\gl_n$ of $\gl_{m|n}$
act on the basis vectors $\ze^{}_{\Tc}$ of $L(\la)^+_{\mu}$
by the same formulas as
the elements of the subalgebra $\gl_n$ of $\gl_{r+n}$
act on the basis vectors $\ze^{}_{\Vc}$ of $\BL(\la')^+_{\mu'}$,
thus completing the proof.
\epf

\bco\label{cor:bact}
Let $\Tc$ be a row-strict tableau of shape
$\Ga_{\la}/\mu$. Then for $p=1,\dots,n$
\beql{akuact}
A_p(u)\ts \ze^{}_{\Tc}=
\frac{(u+l^{\tss\prime}_{r+p,1})\dots (u+l^{\tss\prime}_{r+p,r+p})}
{(u+l^{\tss\prime}_{r+p-1,1})\dots (u+l^{\tss\prime}_{r+p-1,r+p-1})}
\ts \ze^{}_{\Tc}.
\eeq
Moreover, for $p\in\{1,\dots,n-1\}$ and $j\in\{1,\dots,r+p\}$
we have
\beql{cklact}
C_p(-l^{\tss\prime}_{r+p,j})\tss \ze^{}_{\Tc}=\ze^{}_{\Tc^{-}},
\eeq
if $\Tc$ contains an entry $m+p$ in column $j$ and the replacement
of this entry by $m+p+1$ yields a row-strict tableau
$\Tc^{-}$.

For the same values of the parameters $p$ and $j$
we have
\beql{bklact}
B_p(-l^{\tss\prime}_{r+p,j})\tss \ze^{}_{\Tc}=\ze^{}_{\Tc^{+}},
\eeq
if $\Tc$ contains an entry $m+p+1$ in column $j$ and the replacement
of this entry by $m+p$ yields a row-strict tableau
$\Tc^{+}$.
\eco

\bpf
We argue as in the proof
of Proposition~\ref{prop:basskew}.
The rational function $A_p(u)$ coincides with the image
of the series
\ben
{}t^{\ts 1,\dots,\tss p}_{\ts 1,\dots,\tss p}(u)
\ts
t^{\ts 1,\dots,\tss p-1}_{\ts 1,\dots,\tss p-1}(u)^{-1}
\tss(u-h_1-p+1)
\een
under the homomorphism of Proposition~\ref{prop:ymzahom}.
The image of this series (with $h_1$ replaced by $r$)
in the skew representation
$\BL(\la')^+_{\mu'}$ is given by
\ben
A_p(u)=
\big(u+\BE\big)^{1,\dots,\tss r+p}_{1,\dots,\tss r+p}
\cdot
\Big[\big(u+\BE\big)^{1,\dots,r+p-1}_{1,\dots,r+p-1}\Big]^{-1}.
\een
Hence, using again the formulas for the Yangian action in the
Gelfand--Tsetlin basis of $\BL(\la')^+_{\mu'}$
we get the first relation.

The claim involving the operators $C_p(-l^{\tss\prime}_{r+p,j})$
was established in the proof of Proposition~\ref{prop:basskew}.
Similarly,
the rational function $B_p(u)$ coincides with the image
of the series
\ben
{}-t^{\ts 1,\dots,\tss p}_{\ts 1,\dots,\tss p-1,\tss p+1}(u)
\ts
t^{\ts 1,\dots,\tss p+1}_{\ts 1,\dots,\tss p+1}(u)^{-1}
\tss\frac{1}{u-h_1-p}
\een
under the homomorphism of Proposition~\ref{prop:ymzahom}.
The image of this series in
$\BL(\la')^+_{\mu'}$ is given by
\ben
B_p(u)=-
\big(u+\BE\big)^{1,\dots,\tss r+p}_{1,\dots,\tss r+p-1, \tss r+p+1}
\cdot
\Big[\big(u+\BE\big)^{1,\dots,r+p+1}_{1,\dots,r+p+1}\Big]^{-1}.
\een
Let $\Vc$ be the trapezium pattern corresponding to the supertableau
$\Tc$. The formulas for the Yangian action in the
Gelfand--Tsetlin basis of $\BL(\la')^+_{\mu'}$
now give
\ben
B_p(-l^{\tss\prime}_{r+p,j})\tss \ze^{}_{\Vc}=\ze^{}_{\Vc^{+}},
\een
where $\Vc^{+}$ is the pattern obtained from $\Vc$
by replacing the entry $\la^{\tss\prime}_{r+p,j}$
by $\la^{\tss\prime}_{r+p,j}+1$;
see \cite[Sec.~5.4]{m:yc}. It is clear that
$\Vc^{+}$ corresponds to the supertableau
$\Tc^{+}$.
\epf

Corollary~\ref{cor:bact} implies the following identities
for the parameters of the initial tableau $\Tc^{\tss0}$
of shape $\Ga_{\la}/\mu$.

\ble\label{lem:ida}
For each $p=1,\dots,n$ we have
\begin{multline}
\frac{(u+l^{\tss0}_{r+p,1}+p-1)\dots (u+l^{\tss0}_{r+p,r+p}+p-1)}
{(u+l^{\tss0}_{r+p-1,1}+p-1)\dots (u+l^{\tss0}_{r+p-1,r+p-1}+p-1)}\\
{}=\frac{(u+\la_{m+p}+m)(u-r)}{u+m}\ts
\prod_{j=1,\ \la_j-\mu_j\geqslant p}^m\ts \frac{u-\si_j+1}{u-\si_j}.
\non
\end{multline}
\ele

\bpf
Calculate $A_p(u+p-1)\ts\ze_{\Tc^{\tss0}}$ in two different ways
and compare the eigenvalues. First
take $\Tc=\Tc^{\tss0}$ in \eqref{akuact} and replace $u$ by $u+p-1$.
This gives the left hand side of the equality.
On the other hand, since $\ze_{\Tc^{\tss0}}=\ze_{\mu}$ is the
highest vector of the $\Y(\gl_n)$-module $L(\la)^+_{\mu}$,
the eigenvalue of the operator
$Z^{\ts m+1,\dots,\tss m+p}_{\ts m+1,\dots,\tss m+p}(u)$
on $\ze_{\mu}$ can be found from \eqref{eigenzemu}
and \eqref{quaminhv}. Hence, the application of
\eqref{defaku} yields the eigenvalue of $A_p(u+p-1)$
in the form of the right hand side of the equality.
\epf

Using the Gelfand--Tsetlin basis \cite{gt:fdu} of each representation
$L'(\mu)$ of $\gl_m$ with $\mu$
satisfying the assumptions of Corollary~\ref{cor:nonzero},
and the basis $\ze^{}_{\Tc}$
of $L(\la)^+_{\mu}$ formed by the row-strict tableaux $\Tc$
of shape $\Ga_{\la}/\mu$,
we will construct basis vectors of the
vector spaces $L'(\mu)\ot L(\la)^+_{\mu}$ occurring in
\eqref{decomptp}.
More precisely, consider the raising and lowering
operators $s_{ik}$ and $s_{ki}$ for $k=2,\dots,m$ and
$i=1,\dots,k-1$ which are elements of $\U(\gl_m)$
given by the explicit formulas analogous to
\eqref{simae} and \eqref{smaie},

\begin{align}
s_{ik}&=\sum_{i>i_1>\dots>i_p\geqslant 1}
E_{ii_1}E_{i_1i_2}\dots E_{i_{p-1}i_p}E_{i_pk}
(h_i-h_{j_1})\dots (h_i-h_{j_r}),
\non\\
s_{ki}&=\sum_{i<i_1<\dots<i_p<k}
E_{i_1i}E_{i_2i_1}\dots E_{i_pi_{p-1}}E_{ki_p}
(h_i-h_{j_1})\dots (h_i-h_{j_r}),
\non
\end{align}
where $p$ runs over nonnegative integers, $h_i=E_{ii}-i+1$
and $\{j_1,\dots,j_r\}$
is the complementary subset to $\{i_1,\dots,i_p\}$
in the set $\{1,\dots,i-1\}$ or $\{i+1,\dots,k-1\}$, respectively.
Note that the lowering operators $s_{ki}$ and $s_{kj}$
commute for any $i,j\in\{1,\dots,k-1\}$.

Any supertableau $\La$ of shape $\Ga_{\la}$
uniquely determines a partition $\mu=(\mu_1,\dots,\mu_m)$ as the shape
of the subtableau with entries
in $\{1,\dots,m\}$, as well as the corresponding pattern $\Uc$
with entries $\la_{ij}$; see Sec.~\ref{sec:mt}.
Moreover, a row-strict tableau $\Tc$ is found as
the subtableau of $\La$ with entries
in $\{m+1,\dots,m+n\}$. Define the vector of $L(\la)$
corresponding to $\La$ by
\beql{basvect}
\ze^{}_{\La}=\prod_{k=2,\dots,m}^{\longrightarrow}
\Big(s_{k1}^{\la^{}_{k1}-\la^{}_{k-1,1}}\dots
s_{k,k-1}^{\la^{}_{k,k-1}-\la^{}_{k-1,k-1}}\Big)\ts
\ze^{}_{\Tc},
\eeq
where $\ze^{}_{\Tc}$ is defined in \eqref{basvdef}. The following
is implied by the properties of the Gelfand--Tsetlin bases;
see e.g.~\cite{m:gtb} for a proof.

\bpr\label{prop:gt}
The action of the generators of the Lie subalgebra $\gl_m$
of $\gl_{m|n}$ on the vectors $\ze^{}_{\La}$ of $L(\la)$
is given by formulas \eqref{ekk} with $1\leqslant s\leqslant m$
and by formulas \eqref{ekkp}, \eqref{ekkm} with
$1\leqslant s\leqslant m-1$.
\qed
\epr

Due to Propositions~\ref{prop:basskew} and \ref{prop:gt},
in order to determine the action of all elements of $\gl_{m|n}$
on the vectors $\ze^{}_{\La}$ of $L(\la)$,
it will be sufficient to find
explicit expansions of $E_{m,m+1}\ts\ze^{}_{\La}$ and
$E_{m+1,m}\ts\ze^{}_{\La}$ as linear combinations of these vectors.
The following two lemmas are well-known in the
classical case (i.e. for the Lie algebra $\gl_{m+1}$)
and their proofs are not essentially different in the super case;
cf. \cite[Sec.~2.3]{m:gtb} and \cite{z:gz}.
Here we regard the raising and lowering
operators given by \eqref{simae} and \eqref{smaie}
as elements of the universal enveloping algebra
$\U(\gl_{m|n})$.

\ble\label{lem:emmpo}
The following
relation holds in $\U'(\gl_{m|n})$,
\ben
E_{m,m+1}=\sum_{i=1}^{m}s_{m\tss i}\ts z_{i,m+1}
\frac{1}{(h_i-h_1)\ldots\wedge_i\ldots(h_i-h_{m})},
\een
where $s_{mm}=1$.
\qed
\ele

\ble\label{lem:empom}
For any nonnegative integers $k_1,\dots,k_{m-1}$,
in $\U'(\gl_{m|n})$ we have
\ben
\bal
E_{m+1,m}\ts s_{m\tss 1}^{k_1}&\dots s_{m,m-1}^{k_{m-1}}\\
&{}=\sum_{i=1}^{m}
s_{m\tss 1}^{k_1}\dots s_{m\tss i}^{k_i-1}\dots
s_{m,m-1}^{k_{m-1}}\ts z_{m+1,i}
\frac{(h_i-h_1+k_1)\dots(h_i-h_{m-1}+k_{m-1})}
{(h_i-h_1)\ldots\wedge_i\ldots(h_i-h_{m})}.
\eal
\een
\vskip-1.0\baselineskip
\qed
\ele

Observe that $E_{m,m+1}$ commutes with all lowering
operators $s_{ki}$. Therefore, using \eqref{basvect} we get
\beql{emmpozela}
E_{m,m+1}\ts\ze^{}_{\La}=
\prod_{k=2,\dots,m}^{\longrightarrow}
\Big(s_{k1}^{\la^{}_{k1}-\la^{}_{k-1,1}}\dots
s_{k,k-1}^{\la^{}_{k,k-1}-\la^{}_{k-1,k-1}}\Big)\ts
E_{m,m+1}\ts
\ze^{}_{\Tc}.
\eeq
By Lemma~\ref{lem:emmpo}, to find the expansion of
$E_{m,m+1}\ts\ze^{}_{\La}$ we need to calculate
$z_{i,m+1}\ts \ze^{}_{\Tc}$ in terms of
the basis vectors of $L(\la)^+_{\mu+\de_i}$.
Similarly, $E_{m+1,m}$ commutes with the lowering
operators $s_{ki}$ for $k\leqslant m-1$. Hence,
due to Lemma~\ref{lem:empom},
to find the expansion of
$E_{m+1,m}\ts\ze^{}_{\La}$ we need to calculate
$z_{m+1,i}\ts \ze^{}_{\Tc}$ in terms of
the basis vectors of $L(\la)^+_{\mu-\de_i}$.

\ble\label{lem:repbc} For any $p\in\{1,\dots,n-1\}$ we have
in $\Z(\gl_{m|n},\gl_m)${\rm:}
\ben
B_p(u)\tss z_{m+1,i}=z_{m+1,i}\tss B_p(u),\qquad\text{for}
\quad i\geqslant 2,
\een
and
\ben
B_p(u)\tss z_{m+1,1}=z_{m+1,1}\tss B_p(u)\tss
\frac{u-h_1-p}{u-h_1-p+1}.
\een
\ele

\bpf
It suffices to apply \eqref{aux} to permute
$z_{m+1,i}$ with $B_p(u)$, and use the relation
$
z_{m+1,i}\ts(u-h_i-p+1)=(u-h_i-p)\ts z_{m+1,i}
$
for $i\geqslant 2$.
\epf

\ble\label{lem:repcz}
We have the relations in $\Z(\gl_{m|n},\gl_m)${\rm:}
\ben
\bal
z_{1,m+1}\ts C_1(u)&=C_1(u)\ts z_{1,m+1}\ts
\frac{u-h_1+1}{u-h_1-1},\\
z_{i,m+1}\ts C_1(u)&=C_1(u)\ts z_{i,m+1}\ts
\frac{u-h_i+1}{u-h_i},\qquad\text{if}
\quad i\geqslant 2,
\eal
\een
and
\begin{alignat}{2}
z_{i,m+1}\ts C_p(u)&=C_p(u)\ts z_{i,m+1},\qquad&&\text{if}
\quad i,p\geqslant 2,
\non\\
z_{1,m+1}\ts C_p(u)&=C_p(u)\ts z_{1,m+1}\ts
\frac{u-h_1-p+1}{u-h_1-p},\qquad&&\text{if}
\quad p\geqslant 2.
\non
\end{alignat}
\ele

\bpf
All relations follow by the application of \eqref{auxtwo}.
\epf

\bpr\label{prop:zimpe} For any $i\in\{1,\dots,m\}$ we have
\beql{zimck}
z_{i,m+1}\ts\ze^{}_{\Tc}=b_{i,\Tc}\ts\prod_{(p,j)}
\Big(C_p(-l^{\tss\prime}_{r+p,j}-1)
\dots C_p(-l^{\tss 0}_{r+p,j}+1)\tss
C_p(-l^{\tss 0}_{r+p,j})\Big)\ts z_{i,m+1}\ts\ze_{\mu},
\eeq
where the product over the pairs $(p,j)$ is
taken in the same ordering
as in \eqref{basvdef} and
$b_{i,\Tc}$ is a constant given by
\begin{align}
\label{zonempoze}
b_{i,\Tc}&=\prod_{p=1}^{n-1}\ts\prod_{j=1}^{r+p}\ts
\frac{l^{\tss\prime}_{r+p,j}+\si_1+p}{l^{\tss 0}_{r+p,j}+\si_1+p}
\cdot\prod_{j=1}^{r+1}\ts
\frac{l^{\tss\prime}_{r+1,j}+\si_1}{l^{\tss 0}_{r+1,j}+\si_1}
\qquad\text{if}\quad i=1,\\[-1em]
\intertext{and}
\label{zimpoze}
b_{i,\Tc}&=\prod_{j=1}^{r+1}\ts
\frac{l^{\tss\prime}_{r+1,j}+\si_i}{l^{\tss 0}_{r+1,j}+\si_i}
\qquad\text{if}\quad i\geqslant 2.
\end{align}
\epr

\bpf
To verify \eqref{zimpoze}, note that if $i\geqslant 2$, then by
Lemma~\ref{lem:repcz}, $z_{i,m+1}$ is permutable
with the operators of the form $C_k(u)$ with $k\geqslant 2$.
Furthermore, the lemma implies that for $p\geqslant 0$
\ben
\bal
z_{i,m+1}&\ts C_1(u+p-1)
\dots C_1(u+1)\tss
C_1(u)\\
{}&=C_1(u+p-1)
\dots C_1(u+1)\tss
C_1(u)\ts z_{i,m+1}\ts
\frac{u-h_i+p}{u-h_i}
\eal
\een
which yields \eqref{zimpoze}. Relation
\eqref{zonempoze} is verified by a similar calculation with the use
of Lemma~\ref{lem:repcz}.
\epf

Due to Proposition~\ref{prop:zimpe}, the calculation
of $z_{i,m+1}\ts\ze^{}_{\Tc}$ is reduced to expanding
of $z_{i,m+1}\ts\ze_{\mu}$.
If $\la_i-\mu_i=0$ for some $i\in\{1,\dots,m\}$, then
$z_{i,m+1}\ts\ze_{\mu}=0$ by Corollary~\ref{cor:zact}.
On the other hand, if $k=\la_i-\mu_i\geqslant 1$, then
applying \eqref{zmze} with $\mu$ replaced by $\mu+\de_i$
we obtain
\beql{zemuzemp}
\ze_{\mu}=c_{i,\mu}\ts z_{m+k,i}\ts\ze_{\tss\mu+\de_i},
\eeq
where
\ben
c_{i,\mu}=
\prod_{j=1}^{i-1}(-1)^{\la_j-\mu_j}
\frac{\si_i-l_j}{\si_i-\si_j}\ts
\prod_{j=1,\ \la_j-\mu_j\geqslant k}^{i-1}
\frac{\si_i-\si_j}{\si_i-\si_j+1}.
\een
Furthermore, using \eqref{zabhimo}, we derive
\ben
z_{i,m+1}\ts z_{m+k,i}\ts\ze_{\tss\mu+\de_i}
=-Z_{m+k,m+1}(h_i-1)\ts \ze_{\tss\mu+\de_i}=
-Z_{m+k,m+1}(\si_i)\ts \ze_{\tss\mu+\de_i}.
\een
If $k=1$, then by \eqref{eigenzemu}, this vector
equals $\ze_{\tss\mu+\de_i}$ multiplied by a scalar.
Now suppose that $k\geqslant 2$.
In this case relations \eqref{Ess}
imply that
\ben
Z_{m+k,m+1}(\si_i)=[E_{m+k,m+1},Z_{m+1,m+1}(\si_i)],
\een
and so
\beql{zizmp}
z_{i,m+1}\ts z_{m+k,i}\ts\ze_{\tss\mu+\de_i}=
Z_{m+1,m+1}(\si_i)\ts E_{m+k,m+1}\ts\ze_{\tss\mu+\de_i}
\eeq
since $Z_{m+1,m+1}(\si_i)\ts\ze_{\tss\mu+\de_i}=0$ by
\eqref{eigenzemu}.
We will identify $\ze_{\tss\mu+\de_i}$
with the basis vector $\ze^{}_{\Tc^{\tss 0}}$ of
$L(\la)^+_{\mu+\de_i}$ corresponding
to the initial tableau $\Tc^{\tss 0}$ of shape $\Ga_{\la}/(\mu+\de_i)$,
obtained by filling in the boxes of each row by
the consecutive numbers $m+1,m+2,\dots$ from left to right.
We suppose first that $i\geqslant 2$ and
use the corresponding parameters $\la_{r+p,j}^0$
of $\Tc^{\tss 0}$ defined in \eqref{lazero} together with
$l_{r+p,j}^{\tss0}$ defined in \eqref{lzero}.

As we showed in Proposition~\ref{prop:basskew},
the expansion of $E_{m+k,m+1}\ts\ze^{}_{\Tc^{\tss 0}}$
in terms of the basis vectors $\ze^{}_{\Tc}$
can be found by the respective case of the formula \eqref{ekkm}.
Writing
\beql{ecomm}
E_{m+k,m+1}=[E_{m+k,m+k-1},
\dots,[E_{m+3,m+2},E_{m+2,m+1}]\dots]
\eeq
we can see that for $i\geqslant 2$ the expansion of
$E_{m+k,m+1}\ts\ze^{}_{\Tc^{\tss 0}}$ will contain
only those vectors $\ze^{}_{\Tc}$ for which exactly one of
the parameters $\la_{r+1,1}^0,\dots,\la_{r+1,r+1}^0$
is decreased by $1$.
Due to the subsequent application of the operator
$Z_{m+1,m+1}(\si_i)$, there will be only one
of such vectors
occurring in the expansion of \eqref{zizmp}
with a nonzero coefficient.
Indeed, using \eqref{akuact} with $p=1$ we derive that
for any skew tableau $\Tc$ of shape $\Ga_{\la}/(\mu+\de_i)$,
\begin{align}\label{zmmmts}
Z_{m+1,m+1}(u)\ts\ze^{}_{\Tc}&=
\frac{(u+l^{\tss\prime}_{r+1,1})\dots (u+l^{\tss\prime}_{r+1,r+1})}
{(u+l^{\tss\prime}_{r,1})\dots (u+l^{\tss\prime}_{r,r})}\\
\non
{}&\times (u-\si_2)\dots(u-\si_i-1)\dots(u-\si_{m+1})
\ts \ze^{}_{\Tc}.
\end{align}
Set $s=\mu_i$ so that $\mu'_{s+1}=i-1$. If
$\Tc$ occurs in the expansion of $E_{m+k,m+1}\ts\ze^{}_{\Tc^{\tss 0}}$
and $\la'_{r+1,s+2}=\mu'_{s+1}+1$, then
\ben
l^{\tss\prime}_{r+1,s+2}=\la'_{r+1,s+2}-s-1=\mu'_{s+1}-s=i-1-\mu_i
=-\si_i
\een
and $Z_{m+1,m+1}(\si_i)\ts\ze^{}_{\Tc}=0$. Hence,
if $Z_{m+1,m+1}(\si_i)\ts\ze^{}_{\Tc}\ne0$ then
$\la'_{r+1,s+2}=\mu'_{s+1}$.
The betweenness conditions then give
\ben
\la'_{r+1,s+2}=\dots=\la'_{r+k-1,s+k}=\mu'_{s+1},
\een
while for the remaining parameters we have
$\la'_{r+p,j}=\la^0_{r+p,j}$. This determines a unique tableau
which we denote by $\Tc^{+}$. Thus,
\beql{ziztp}
z_{i,m+1}\ts z_{m+k,i}\ts\ze_{\tss\mu+\de_i}
=d_{i,\mu}\ts\ze^{}_{\Tc^{+}},
\eeq
for a nonzero constant $d_{i,\mu}$ and any $k\geqslant 1$.
To calculate its value, note that by \eqref{ekkm} and \eqref{ecomm}
the coefficient of the basis vector $\ze^{}_{\Tc^{+}}$
in the expansion of $E_{m+k,m+1}\ts\ze_{\tss\mu+\de_i}$
coincides with the coefficient of this vector
in the expansion of
\ben
(-1)^k\ts
E_{m+2,m+1}\ts E_{m+3,m+2}\dots E_{m+k,m+k-1}\ts\ze_{\tss\mu+\de_i}
\een
and it is found by the formula
\beql{coeflow}
(-1)^k\ts\prod_{p=1}^{k-1}\ts
\frac{(\si_i+l^{\tss0}_{r+p-1,1}+p-1)\dots
(\si_i+l^{\tss0}_{r+p-1,\tss r+p-1}+p-1)}
{(\si_i+l^{\tss0}_{r+p,1}+p-1)\ldots
\wedge_{s+p+1}\ldots(\si_i+l^{\tss0}_{r+p,\tss r+p}+p-1)}.
\eeq
We now use Lemma~\ref{lem:ida} (applied to $\Ga_{\la}/(\mu+\de_i)$
instead of $\Ga_{\la}/\mu$)
to write this coefficient
in a different form. Divide both sides of the identity
of the lemma by $u-\si_i$ and set $u=\si_i$.
We get the identity
\begin{multline}
\frac{(\si_i+l^{\tss0}_{r+p-1,1}+p-1)\dots
(\si_i+l^{\tss0}_{r+p-1,\tss r+p-1}+p-1)}
{(\si_i+l^{\tss0}_{r+p,1}+p-1)\ldots
\wedge_{s+p+1}\ldots(\si_i+l^{\tss0}_{r+p,\tss r+p}+p-1)}\\
{}=-\frac{\si_i+m}{(\si_i-r)(\si_i+\la_{m+p}+m)}\ts
\prod_{j=1,\ j\ne i,\ \la_j-\mu_j\geqslant p}^m\ts
\frac{\si_i-\si_j}{\si_i-\si_j+1}.
\non
\end{multline}
Taking the product over $p$, we find that
the coefficient \eqref{coeflow} equals
\ben
{}-\Big(\frac{\si_i+m}{\si_i-r}\Big)^{k-1}\ts
\prod_{p=1}^{k-1}\frac{1}{\si_i+\la_{m+p}+m}
\ts\prod_{j=1,\ j\ne i}^m\ts
\Big(\frac{\si_i-\si_j}
{\si_i-\si_j+1}\Big)^{\min\{\la_j-\mu_j,\ts k-1\}}.
\een
Hence, applying \eqref{zmmmts} for $u=\si_i$ and $\Tc=\Tc^{+}$,
and combining this with \eqref{zemuzemp},
we conclude that in the case under consideration,
\beql{zimglow}
z_{i,m+1}\ts \ze_{\tss\mu}
=g_{i,\mu}\ts\ze^{}_{\Tc^{+}},
\eeq
with
\begin{multline}
g_{i,\mu}
=\prod_{j=1}^{i-1}(-1)^{\la_j-\mu_j}
\frac{\si_i-l_j}{\si_i-\si_j}\ts
\prod_{j=1,\ \la_j-\mu_j\geqslant k}^{i-1}
\frac{\si_i-\si_j}{\si_i-\si_j+1}\\[1em]
{}\times\Big(\frac{\si_i+m}{\si_i-r}\Big)^{k-1}\ts
\prod_{p=1}^{k-1}\frac{1}{\si_i+\la_{m+p}+m}
\ts\prod_{j=1,\ j\ne i}^m\ts
\Big(\frac{\si_i-\si_j}
{\si_i-\si_j+1}\Big)^{\min\{\la_j-\mu_j,\ts k-1\}}\\[1em]
{}\times\frac{(\si_i+l^{\tss0}_{r+1,1})\dots
(\si_i+l^{\tss0}_{r+1,\tss r+1})}
{(\si_i+l^{\tss0}_{r,1})\ldots
\wedge_{s+1}\ldots(\si_i+l^{\tss0}_{r,\tss r})}
\ts (\si_i-\si_2)\ldots\wedge_{i}\ldots (\si_i-\si_{m+1}),
\non
\end{multline}
where the parameters $l^{\tss0}_{r+p,\tss j}$ are now
associated with the initial tableau $\Tc^{\tss0}$
of shape $\Ga_{\la}/\mu$ via \eqref{lazero} and \eqref{lzero}.
The same formula for $g_{i,\mu}$ clearly remains
valid in the case $k=1$ as well.

To extend the above calculation to the case $i=1$
we need to take into account the fact that the parameter $r=\mu_1$
changes to $r+1$ for the partition $\mu+\de_1$.
This time \eqref{zmmmts} is replaced by the relation
\beql{zmmmtspo}
Z_{m+1,m+1}(u)\ts\ze^{}_{\Tc}=
\frac{(u+l^{\tss\prime}_{r+2,1})\dots (u+l^{\tss\prime}_{r+2,r+2})}
{(u+l^{\tss\prime}_{r+1,1})\dots (u+l^{\tss\prime}_{r+1,r+1})}\ts
(u-\si_2)\dots(u-\si_{m+1})
\ts \ze^{}_{\Tc},
\eeq
where $\Tc$ is any skew
tableau of shape $\Ga_{\la}/(\mu+\de_1)$.
We will use the parameters $\la^0_{r+p,j}$ and
$l^{\tss0}_{r+p,j}$ associated
with the initial tableau of shape $\Ga_{\la}/\mu$
by \eqref{lazero} and \eqref{lzero}.
Relation
\eqref{ziztp} holds for $i=1$ as well, where
the parameters of the tableau $\Tc^{+}$ of shape
$\Ga_{\la}/(\mu+\de_1)$ are given by
\beql{parta}
\la'_{r+p+1,j}=\la^0_{r+p,j}\qquad\text{for}\quad
p=0,\dots,n \quad\text{and}\quad
j=1,\dots,r+p,
\eeq
while $\la'_{r+1,r+1}=1$ and $\la'_{r+p,r+p}=0$
for $p\geqslant 2$.
The coefficient of the basis vector $\ze^{}_{\Tc^{+}}$
in the expansion of $E_{m+k,m+1}\ts\ze_{\tss\mu+\de_1}$
coincides with the coefficient of this vector
in the expansion of
\ben
(-1)^k\ts
E_{m+2,m+1}\ts E_{m+3,m+2}\dots E_{m+k,m+k-1}\ts\ze_{\tss\mu+\de_1}
\een
and it is found by the formula
\beql{coeflowone}
(-1)^k\ts\prod_{p=1}^{k-1}\ts
\frac{(r+l^{\tss0}_{r+p-1,1}+p-1)\dots
(r+l^{\tss0}_{r+p-1,\tss r+p-1}+p-1)}
{(r+l^{\tss0}_{r+p,1}+p-1)\dots
(r+l^{\tss0}_{r+p,\tss r+p}+p-1)}.
\eeq
Lemma~\ref{lem:ida} now gives
\begin{multline}
\frac{(r+l^{\tss0}_{r+p-1,1}+p-1)\dots
(r+l^{\tss0}_{r+p-1,\tss r+p-1}+p-1)}
{(r+l^{\tss0}_{r+p,1}+p-1)\dots
(r+l^{\tss0}_{r+p,\tss r+p}+p-1)}\\
{}=\frac{r+m}{r+\la_{m+p}+m}\ts
\prod_{j=2,\ \la_j-\mu_j\geqslant p}^m\ts
\frac{r-\si_j}{r-\si_j+1},
\non
\end{multline}
so that
taking the product over $p$, we find that
the coefficient \eqref{coeflowone} equals
\ben
(-1)^k\ts(r+m)^{k-1}\ts
\prod_{p=1}^{k-1}\frac{1}{r+\la_{m+p}+m}
\ts\prod_{j=2}^m\ts
\Big(\frac{r-\si_j}
{r-\si_j+1}\Big)^{\min\{\la_j-\mu_j,\ts k-1\}}.
\een
Therefore, applying \eqref{zmmmtspo} for $u=r$ and $\Tc=\Tc^{+}$,
we conclude that \eqref{zimglow}
holds for $i=1$
with
\begin{multline}
g_{1,\mu}
=(-1)^{k-1}\ts(r+m)^{k-1}\ts
\prod_{p=1}^{k-1}\frac{1}{r+\la_{m+p}+m}
\ts\prod_{j=2}^m\ts
\Big(\frac{r-\si_j}
{r-\si_j+1}\Big)^{\min\{\la_j-\mu_j,\ts k-1\}}\\[1em]
{}\times\frac{(r+l^{\tss0}_{r+1,1})\dots
(r+l^{\tss0}_{r+1,\tss r+1})}
{(r+l^{\tss0}_{r,1})\dots(r+l^{\tss0}_{r,\tss r})}
\ts (r-\si_2)\dots (r-\si_{m+1}),
\non
\end{multline}
which is valid for $k\geqslant 1$.

Combining the above calculation with Proposition~\ref{prop:zimpe}
we come to the following.

\bpr\label{prop:zimmpa}
Suppose that for some $i\in\{1,\dots,m\}$ the following
condition holds{\rm:}
$\mu+\de_i$ is a partition,
and a row-strict tableau $\Tc$
of shape $\Ga_{\la}/\mu$ contains
the entry $m+1$ in the box $(i,\mu_i+1)$. Then
\beql{zimponeta}
z_{i,m+1}\ts \ze^{}_{\Tc}=b_{i,\Tc}
\ts g_{i,\mu}\ts \ze^{}_{\Tc^{+}_i},
\eeq
where $\Tc^{+}_i$ is the tableau obtained from $\Tc$
by removing the entry $m+1$ from the box $(i,\mu_i+1)$, and
the coefficients
$b_{i,\Tc}$ and $g_{i,\mu}$ are nonzero.
Moreover, if the condition does not hold, then
$z_{i,m+1}\ts \ze^{}_{\Tc}=0$.
\epr

\bpf
If the condition holds, then the claim follows from
the formulas \eqref{cklact} and
the definition \eqref{basvdef} of the vectors $\ze^{}_{\Tc}$
together with \eqref{zimck} and \eqref{zimglow}.
Now suppose that
$\mu+\de_i$ is not a partition. If the vector
$\xi=z_{i,m+1}\ts \ze^{}_{\Tc}$ were nonzero, its cyclic span
$\U(\gl_m)\ts\xi$ would be a highest weight $\gl_m$-module
with the highest weight $\mu+\de_i$. Since
$\mu+\de_i$ is not a partition, this module must be
infinite-dimensional, a contradiction.

Furthermore, if $\mu+\de_i$ is a partition, then
the condition will be violated if the box $(i,\mu_i+1)$
is outside the diagram $\Ga_{\la}/\mu$. In this case
$\la_i=\mu_i$ and the claim follows from
Corollary~\ref{cor:zact}, as $z_{i,m+1}\ts\ze_{\mu}=0$.
Finally, let the diagram $\Ga_{\la}/\mu$ contain the box $(i,\mu_i+1)$
with the entry of $\Tc$ in this box greater than $m+1$.
Set $j=\mu_i+1$. Then all entries in column $j$ of the tableau
$\Tc$ should also exceed $m+1$ and so the parameter
$l^{\tss\prime}_{r+1,j}$ of $\Tc$ equals
$\la'_{r+1,j}-j+1=-\mu_i+i-1=-\si_i$. Hence, $b_{i,\Tc}=0$
and the claim follows from Proposition~\ref{prop:zimpe}.
\epf

Proposition~\ref{prop:zimmpa} together with Lemma~\ref{lem:emmpo}
and relation \eqref{emmpozela} provide
explicit formulas for the coefficients in the expansion
of $E_{m,m+1}\ts\ze^{}_{\La}$ as a linear
combination of the basis vectors.
We will express these coefficients in terms of the parameters
of $\La$ in Theorem~\ref{thm:baslla} below.

We now turn to the calculation
of $z_{m+1,i}\ts \ze^{}_{\Tc}$ for an arbitrary row-strict
tableau $\Tc$ of shape $\Ga_{\la}/\mu$.
Suppose that the condition of Proposition~\ref{prop:zimmpa}
holds.
Apply $z_{m+1,i}$ to both sides of \eqref{zimponeta}.
Using \eqref{zabhi}, we get
\beql{fitz}
b_{i,\Tc}\ts g_{i,\mu}\ts z_{m+1,i}\ts\ze^{}_{\Tc^{+}_i}=
z_{m+1,i}\ts z_{i,m+1}\ts \ze^{}_{\Tc}
=Z_{m+1,m+1}(h_i)\ts \ze^{}_{\Tc}=Z_{m+1,m+1}(\si_i)\ts \ze^{}_{\Tc}
\eeq
which equals $\ze^{}_{\Tc}$ multiplied by a scalar
found for $i=1$ and $i\geqslant 2$
from \eqref{zmmmtspo} and \eqref{zmmmts}, respectively.
This allows us to calculate $z_{m+1,i}\ts\ze^{}_{\Tc^{+}_i}$.

\bpr\label{prop:lowact}
Suppose that for some $i\in\{1,\dots,m\}$ the following
condition holds{\rm:}
$\mu-\de_i$ is a partition
and a row-strict tableau $\Tc$
of shape $\Ga_{\la}/\mu$ does not contain
the entry $m+1$ in the box $(i,\mu_i+1)$. Then
\beql{zmponei}
z_{m+1,i}\ts \ze^{}_{\Tc}=f_{i,\Tc}\ts \ze^{}_{\Tc^{-}_i},
\eeq
where $\Tc^{-}_i$ is the tableau obtained from $\Tc$
by adding the entry $m+1$ in the box $(i,\mu_i)$, and
$f_{i,\Tc}$ is a constant.
If the condition does not hold, then $z_{m+1,i}\ts \ze^{}_{\Tc}=0$.
\epr

\bpf
The first part of the proposition follows from
\eqref{fitz}. Now suppose that the condition is violated
and $\mu-\de_i$ is not a partition. If $i<m$ then the claim
follows by the same argument as in the proof of
Proposition~\ref{prop:zimmpa} by considering
the $\U(\gl_m)$-cyclic span of
the vector $z_{m+1,i}\ts \ze^{}_{\Tc}$.
If $i=m$ then we must have $\mu_m=0$. We will show that
\beql{zmzeta}
z_{m+1,m}\ts\ze^{}_{\Tc}=0
\eeq
for all row-strict tableaux $\Tc$ of skew shapes $\Ga_{\la}/\mu$
with $\mu$ running over the partitions with $\mu_m=0$.
Using \eqref{akuact} with $p=1$ we derive that
\beql{zcazza}
Z_{m+1,m+1}(u)\ts\ze^{}_{\Tc}=
\frac{(u+l^{\tss\prime}_{r+1,1})\dots (u+l^{\tss\prime}_{r+1,r+1})}
{(u+l^{\tss\prime}_{r,1})\dots (u+l^{\tss\prime}_{r,r})}
\cdot (u-\si_2)\dots(u-\si_{m+1})
\ts \ze^{}_{\Tc}.
\eeq
Therefore, by \eqref{zabhimo} we have
\beql{zmmhtz}
z_{m,m+1}\ts z_{m+1,m}\ts\ze^{}_{\Tc}=
-Z_{m+1,m+1}(h_m-1)\ts \ze^{}_{\Tc}
=-Z_{m+1,m+1}(-m)\ts \ze^{}_{\Tc}=0,
\eeq
where we used the assumption $\mu_m=0$ to conclude that
all parameters $l^{\tss\prime}_{rj}$ in \eqref{zcazza}
do not exceed $m-1$, while $\si_{m+1}=-m$. We proceed
by induction on the weights $\om(\Tc)$
of the vectors $\ze^{}_{\Tc}$. The base of induction
is the case $\ze^{}_{\Tc}=\ze_{\mu}$
with $\mu=(\la_1,\dots,\la_{m-1},0)$ so that
\ben
\om(\Tc)=(\la_1,\dots,\la_{m-1},0\ts|\ts
\la_{m+1}+1,\dots,\la_{m+k}+1,\la_{m+k+1},\dots,\la_{m+n}),
\een
where $k=\la_m$. In this case
$
z_{m+1,m}\ts \ze_{\mu}=z_{m+1,m}\ts z_{m+k,m}\dots z_{m+1,m}\ts\zeta=0
$
since $z_{m+1,m}^2=0$ and
$z_{m+1,m}\ts z_{m+p,m}=-z_{m+1,m}\ts z_{m+p,m}$ by \eqref{siasja}.
Now, given an arbitrary $\Tc$ we will show that
the vector $z_{m+1,m}\ts\ze^{}_{\Tc}$ is annihilated
by all operators $E_{a,a+1}$ with $a=1,\dots,m+n-1$.
Indeed, this is clear for $a=1,\dots,m-1$.
Furthermore, each operator $E_{a,a+1}$ with $a>m$ commutes
with $z_{m+1,m}$ so that
\ben
E_{a,a+1}\ts z_{m+1,m}\ts\ze^{}_{\Tc}=
z_{m+1,m}\ts E_{a,a+1}\ts \ze^{}_{\Tc}
\een
which is zero by the induction hypothesis,
since $E_{a,a+1}\ts \ze^{}_{\Tc}$ is a linear
combination of the vectors $\ze^{}_{\Tc^{\tss\prime}}$
with the weights $\om(\Tc^{\tss\prime})$ exceeding
$\om(\Tc)$.
To calculate $E_{m,m+1}\ts z_{m+1,m}\ts\ze^{}_{\Tc}$
use Lemma~\ref{lem:emmpo}. Using \eqref{ziazbj} we find that
for $i<m$
\ben
z_{i,m+1}\ts z_{m+1,m}\ts\ze^{}_{\Tc}
=-z_{m+1,m}\ts z_{i,m+1}\ts\ze^{}_{\Tc}
\een
which is zero by Proposition~\ref{prop:zimmpa} and
the induction hypothesis. Finally, the case $i=m$
is taken care of by \eqref{zmmhtz}. Thus, if a vector
$z_{m+1,m}\ts\ze^{}_{\Tc}$ were nonzero, it would generate
a proper nonzero submodule of $L(\la)$. This is a contradiction
since $L(\la)$ is irreducible. This proves \eqref{zmzeta}.

It remains to show that $z_{m+1,i}\ts \ze^{}_{\Tc}=0$
in the case where $\mu-\de_i$ is a partition and
a row-strict tableau $\Tc$
of shape $\Ga_{\la}/\mu$ contains
the entry $m+1$ in the box $(i,\mu_i+1)$.
By the first part of the proposition, the vector $\ze^{}_{\Tc}$
can be written as
$\ze^{}_{\Tc}=c\ts z_{m+1,i}\ts \ze^{}_{\overline\Tc}$,
where $c$ is a constant and the tableau $\overline\Tc$ is obtained
from $\Tc$ by removing the entry $m+1$ in the box $(i,\mu_1+1)$.
Then $z_{m+1,i}\ts\ze^{}_{\Tc}=0$ since $z_{m+1,i}^2=0$.
\epf

We will now use
Propositions~\ref{prop:zimmpa}
and \ref{prop:lowact} to produce a basis of $L(\la)$
and to prove the main theorem.
Recall the vectors $\ze^{}_{\La}$ of $L(\la)$ constructed
in \eqref{basvect}. As before, to each supertableau
$\La$ of shape $\Ga_{\la}$
we associate a partition $\mu=(\mu_1,\dots,\mu_m)$ as the shape
of the subtableau with entries
in $\{1,\dots,m\}$. The parameters $l_{r+p,j}^{\tss0}$ of the
initial tableau $\Tc^{\tss 0}$ of shape
$\Ga_{\la}/\mu$ are defined in \eqref{lzero}, where
$r=r(\La)$ is defined by $r=\la_{m1}=\mu_1$.

\bth\label{thm:baslla}
The vectors $\ze^{}_{\La}$ parameterized by all
supertableaux $\La$ of shape $\Ga_{\la}$
form a basis of the representation $L(\la)$ of $\gl_{m|n}$.
Moreover, the action of the generators of $\gl_{m|n}$
is given by \eqref{ekk}--\eqref{ekkm} together with
the formulas
\begin{align}
\label{ekkpodd}
E_{m,m+1}\ts \ze^{}_{\Lambda}
&=\sum_{\La'} c_{\La\La'}\ts
\ze^{}_{\La'},\\
\label{ekkmodd}
E_{m+1,m}\ts \ze^{}_{\Lambda}
&=\sum_{\La'} d_{\La\La'}\ts
\ze^{}_{\La'},
\end{align}
where the sums are taken
over supertableaux $\La'$ obtained from $\La$
respectively by replacing
an entry $m+1$ by $m$ and by replacing
an entry $m$ by $m+1$. The coefficients are found by the formulas
\ben
\bal
c_{\La\La'}&=\frac{(l^{}_{mi}+l^{\ts\prime}_{r+1,1})
\dots (l^{}_{mi}+l^{\ts\prime}_{r+1,r+1})}
{(l^{}_{mi}+l^{\tss0}_{r,1})
\ldots\wedge^{}_{\la_{mi}+1}\ldots (l^{}_{mi}+l^{\tss0}_{r,r})}
\ts\Bigg(\frac{l^{}_{mi}+m}
{l^{}_{mi}-l^{}_{m1}}\Bigg)^{k}\\[1em]
{}&\times{}
\prod_{p=1}^{k-1}\frac{1}{l_{mi}+\la_{m+p}+m}
\ts\prod_{j=1}^{i-1}
(-1)^{\la_j-\la_{mj}}\ts\frac{l_{mi}-l_j}{l_{mi}-l_{mj}}
\\[1em]
{}&\times{}\prod_{j=1,\ \la_{j}-\la_{mj}\geqslant k}^{i-1}
\frac{l_{mi}-l_{mj}}{l_{mi}-l_{mj}+1}
\ts\prod_{j=1,\ j\ne i}^{m}\Bigg(\frac{l^{}_{mi}-l^{}_{mj}}
{l^{}_{mi}-l^{}_{mj}+1}\ts
\Bigg)^{\min\{\la_j-\la_{mj},\ts k-1\}}
\eal
\een
and
\ben
\bal
d_{\La\La'}&=\frac{(l^{}_{mi}-l^{}_{m-1,1})
\dots (l^{}_{mi}-l^{}_{m-1,\tss m-1})}
{(l^{}_{mi}-l^{}_{m1})\ldots\wedge_i\ldots(l^{}_{mi}-l^{}_{mm})}
\ts\Bigg(\frac{l^{}_{mi}-l^{}_{m1}-1}
{l^{}_{mi}+m-1}\Bigg)^{k}\\[1em]
{}&\times{}
\prod_{p=1}^{k}(l_{mi}+\la_{m+p}+m-1)
\ts\prod_{j=1}^{i-1}
(-1)^{\la_j-\la_{mj}}\ts\frac{l_{mi}-l_{mj}-1}{l_{mi}-l_j-1}
\\[1em]
{}&\times{}\prod_{j=1,\ \la_{j}-\la_{mj}> k}^{i-1}
\frac{l_{mi}-l_{mj}}{l_{mi}-l_{mj}-1}
\ts\prod_{j=1,\ j\ne i}^{m}\Bigg(\frac{l^{}_{mi}-l^{}_{mj}}
{l^{}_{mi}-l^{}_{mj}-1}\ts
\Bigg)^{\min\{\la_j-\la_{mj},\ts k\}},
\eal
\een
where the
replacement occurs in row
$i\geqslant 2$ and $k=\la_i-\la_{mi}${\rm;}
\ben
\bal
c_{\La\La'}&=(-1)^{k-1}(l_{m1}+m)^{k}
\ts\frac{(l^{}_{m1}+l^{\ts\prime}_{r+1,1})
\dots (l^{}_{m1}+l^{\ts\prime}_{r+1,r+1})}
{(l^{}_{m1}+l^{\tss0}_{r,1})
\dots (l^{}_{m1}+l^{\tss0}_{r,r})}\ts
\prod_{p=1}^{k-1}\frac{1}
{l^{}_{m1}+\la_{m+p}+m}\\[1em]
{}&\times{}\prod_{p=1}^{n-1}\ts\prod_{j=1}^{r+p}\ts
\frac{l^{}_{m1}+l^{\ts\prime}_{r+p,j}+p}
{l^{}_{m1}+l^{\tss0}_{r+p,j}+p}
\ts\prod_{j=2}^{m}\Bigg(\frac{l^{}_{m1}-l^{}_{mj}}
{l^{}_{m1}-l^{}_{mj}+1}\ts
\Bigg)^{\min\{\la_j-\la_{mj},\ts k-1\}}
\eal
\een
and
\ben
\bal
d_{\La\La'}&=\frac{(-1)^k}{(l_{m1}+m-1)^{k}}
\ts\frac{(l^{}_{m1}-l^{}_{m-1,1})
\dots (l^{}_{m1}-l^{}_{m-1,\tss m-1})}
{(l^{}_{m1}-l^{}_{m2})\dots(l^{}_{m1}-l^{}_{mm})}
\\[1em]
{}&\times{}\prod_{p=1}^{k}
(l^{}_{m1}+\la_{m+p}+m-1)
\ts\prod_{j=2}^{m}\Bigg(\frac{l^{}_{m1}-l^{}_{mj}}
{l^{}_{m1}-l^{}_{mj}-1}\ts
\Bigg)^{\min\{\la_j-\la_{mj},\ts k\}}\\[1em]
{}&\times{}\prod_{p=1}^{n-1}\prod_{j=1}^{r+p-1}\ts
\frac{l^{}_{m1}+l^{\tss0}_{r+p,j}+p-1}
{l^{}_{m1}+l^{\ts\prime}_{r+p,j}+p-1},
\eal
\een
where the replacement occurs in row $1$ and $k=\la_1-\la_{m1}$.
\eth

\bpf
Consider the subspace $K$ of $L(\la)$,
spanned by all vectors $\ze^{}_{\La}$.
By Propositions~\ref{prop:basskew} and \ref{prop:gt}
the subspace $K$ is invariant with respect to the action
of the subalgebras $\gl_m$ and $\gl_n$ of $\gl_{m|n}$.
Moreover, Lemmas~\ref{lem:emmpo} and \ref{lem:empom}
together with Propositions~\ref{prop:zimmpa}
and \ref{prop:lowact} imply that $K$ is also invariant
with respect to the action of the elements $E_{m,m+1}$ and
$E_{m+1,m}$. Hence, $K$ is a $\gl_{m|n}$-submodule of $L(\la)$.
Since $K$ contains the highest vector $\ze$, and $L(\la)$
is irreducible we can conclude that $K=L(\la)$.

Furthermore, for each partition $\mu=(\mu_1,\dots,\mu_m)$
such that $0\leqslant\la_i-\mu_i\leqslant n$ for all $i$,
the vectors $\ze_{\Tc}$ parameterized by skew tableaux $\Tc$
of shape $\Ga_{\la}/\mu$ form a basis of the vector space
$L(\la)^+_{\mu}$. For any fixed $\Tc$ the vectors
of the form \eqref{basvect} parameterized by
the column-strict $\mu$-tableaux with entries
in $\{1,\dots,m\}$ form a basis of the $\gl_m$-module
$L'(\mu)$. This implies that the vectors $\ze^{}_{\La}$
are linearly independent and hence form a basis of $L(\la)$.

We will now calculate the expansions $E_{m,m+1}\ts\ze^{}_{\La}$
and $E_{m+1,m}\ts\ze^{}_{\La}$ as linear combinations
of the basis vectors.
Taking into account the denominator in
the relation of Lemma~\ref{lem:emmpo} and using
\eqref{zimglow} for $i\geqslant 2$ we get
\begin{align}
\frac{g_{i,\mu}}{(\si_i-\si_1)\ldots
\wedge_i\ldots(\si_i-\si_m)}
={}&{}\prod_{j=1}^{i-1}(-1)^{\la_j-\mu_j}
\frac{\si_i-l_j}{\si_i-\si_j}\ts
\prod_{j=1,\ \la_j-\mu_j\geqslant k}^{i-1}
\frac{\si_i-\si_j}{\si_i-\si_j+1}
\non\\[1em]
{}\times\Big(\frac{\si_i+m}{\si_i-r}\Big)^{k}\ts
{}&{}\prod_{p=1}^{k-1}\frac{1}{\si_i+\la_{m+p}+m}
\ts\prod_{j=1,\ j\ne i}^m\ts
\Big(\frac{\si_i-\si_j}
{\si_i-\si_j+1}\Big)^{\min\{\la_j-\mu_j,\ts k-1\}}
\non\\[1em]
{}&{}\times\frac{(\si_i+l^{\tss0}_{r+1,1})\dots
(\si_i+l^{\tss0}_{r+1,\tss r+1})}
{(\si_i+l^{\tss0}_{r,1})\ldots
\wedge_{s+1}\ldots(\si_i+l^{\tss0}_{r,\tss r})},
\non
\end{align}
where $s=\mu_i$. The
calculation is completed by
applying Proposition~\ref{prop:zimpe} with formula
\eqref{zimpoze}. Similarly, using the case $i=1$
of \eqref{zimglow}, we find that
\begin{multline}
\frac{g_{1,\mu}}{(r-\si_2)\dots(r-\si_m)}
=(-1)^{k-1}\ts(r+m)^{k}\ts
\prod_{p=1}^{k-1}\frac{1}{r+\la_{m+p}+m}\\[1em]
{}\times\ts\frac{(r+l^{\tss0}_{r+1,1})\dots
(r+l^{\tss0}_{r+1,\tss r+1})}
{(r+l^{\tss0}_{r,1})\dots(r+l^{\tss0}_{r,\tss r})}
\ts\prod_{j=2}^m
\Big(\frac{r-\si_j}
{r-\si_j+1}\Big)^{\min\{\la_j-\mu_j,\ts k-1\}}.
\non
\end{multline}
The calculation is completed by
applying Proposition~\ref{prop:zimpe}
with \eqref{zonempoze} thus proving the formulas for
$c_{\La\La'}$ in \eqref{ekkpodd}.

To prove the formula for the coefficients $d_{\La\La'}$
in \eqref{ekkmodd},
apply $z_{m+1,i}$ with $i\geqslant 2$ to both sides of \eqref{zimglow}
and use \eqref{zabhi} to get
\ben
g_{i,\mu}\ts z_{m+1,i}\ts\ze^{}_{\Tc^{+}}=
z_{m+1,i}\ts z_{i,m+1}\ts \ze_{\tss\mu}
=Z_{m+1,m+1}(h_i)\ts \ze_{\tss\mu}
=Z_{m+1,m+1}(\si_i)\ts \ze_{\tss\mu}
\een
so that
\ben
z_{m+1,i}\ts\ze^{}_{\Tc^{+}}=g^{-1}_{i,\mu}\ts
Z_{m+1,m+1}(\si_i)\ts \ze_{\tss\mu}.
\een
Then
\eqref{zcazza} implies
\ben
Z_{m+1,m+1}(\si_i)\ts \ze_{\tss\mu}=
\frac{(\si_i+l^{\tss0}_{r+1,1})\dots
(\si_i+l^{\tss0}_{r+1,r+1})}
{(\si_i+l^{\tss0}_{r,1})\ldots\wedge_{s+1}\ldots
(\si_i+l^{\tss0}_{r,r})}
\cdot (\si_i-\si_2)\ldots\wedge_i\ldots(\si_i-\si_{m+1})
\ts \ze_{\tss\mu},
\een
where we identify $\ze_{\mu}$ with the vector
$\ze^{}_{\Tc^{\tss0}}$ associated with the initial tableau
$\Tc^{\tss0}$ of the skew
shape $\Ga_{\la}/\mu$. Using the formula for $g_{i,\mu}$ we get
\ben
\bal
z_{m+1,i}\ts\ze^{}_{\Tc^{+}}&=
\Big(\frac{\si_i-r}{\si_i+m}\Big)^{k-1}
\ts\ts\prod_{p=1}^{k-1}(\si_i+\la_{m+p}+m)
\ts\prod_{j=1,\ j\ne i}^m\ts
\Big(\frac{\si_i-\si_j+1}
{\si_i-\si_j}\Big)^{\min\{\la_j-\mu_j,\ts k-1\}}\\
{}&\times{}\ts\prod_{j=1}^{i-1}
(-1)^{\la_j-\mu_j}\ts\frac{\si_i-\si_j}{\si_i-l_j}\ts
\prod_{j=1,\ \la_j-\mu_j\geqslant k}^{i-1}
\frac{\si_i-\si_j+1}{\si_i-\si_j}\ts\ze_{\mu}.
\eal
\een
Replace $\mu$ by $\mu-\de_i$
so that $k$ will be replaced by $k+1$ and
$\Tc^{+}$ will be associated with $\mu$.
Thus, for $k\geqslant 0$ we have
\ben
\bal
z_{m+1,i}\ts\ze^{}_{\Tc^{+}}&=
\Big(\frac{\si_i-r-1}{\si_i+m-1}\Big)^{k}
\ts\ts\prod_{p=1}^{k}(\si_i+\la_{m+p}+m-1)
\ts\prod_{j=1,\ j\ne i}^m\ts
\Big(\frac{\si_i-\si_j}
{\si_i-\si_j-1}\Big)^{\min\{\la_j-\mu_j,\ts k\}}\\
{}&\times{}\ts\prod_{j=1}^{i-1}
(-1)^{\la_j-\mu_j}\ts\frac{\si_i-\si_j-1}{\si_i-l_j-1}\ts
\prod_{j=1,\ \la_j-\mu_j> k}^{i-1}
\frac{\si_i-\si_j}{\si_i-\si_j-1}\ts\ze_{\mu-\de_i}.
\eal
\een
Now we use \eqref{bklact}
and the first relation of Lemma~\ref{lem:repbc}
to calculate $z_{m+1,i}\ts\ze^{}_{\Tc}$ for an arbitrary
skew tableau of shape $\Ga_{\la}/\mu$.
Applying appropriate operators
of the form $B_p(-l^{\tss\prime}_{r+p,j})$
to both sides of \eqref{zmponei}, we conclude that
the coefficient $f_{i,\Tc}$ coincides
with $f_{i,\Tc^+}$, where the tableau $\Tc^+$
is associated with $\mu$ as above. Together with
Lemma~\ref{lem:empom} this proves the formula
for the coefficient $d_{\La\La'}$ in \eqref{ekkmodd},
for $i\geqslant 2$.

Finally, consider the case $i=1$ of the formula \eqref{ekkmodd}.
Apply $z_{m+1,1}$ to both sides of \eqref{zimglow} with $i=1$
and use \eqref{zabhi} to get
\ben
g_{1,\mu}\ts z_{m+1,1}\ts\ze^{}_{\Tc^{+}}=
z_{m+1,1}\ts z_{1,m+1}\ts \ze_{\tss\mu}
=Z_{m+1,m+1}(h_1)\ts \ze_{\tss\mu}
=Z_{m+1,m+1}(r)\ts \ze_{\tss\mu}
\een
so that
\ben
z_{m+1,1}\ts\ze^{}_{\Tc^{+}}=g^{-1}_{1,\mu}\ts
Z_{m+1,m+1}(r)\ts \ze_{\tss\mu}.
\een
Now \eqref{zcazza} gives
\ben
Z_{m+1,m+1}(r)\ts \ze_{\tss\mu}=
\frac{(r+l^{\tss0}_{r+1,1})\dots
(r+l^{\tss0}_{r+1,r+1})}
{(r+l^{\tss0}_{r,1})\dots
(r+l^{\tss0}_{r,r})}
\cdot (r-\si_2)\dots(r-\si_{m+1})
\ts \ze_{\tss\mu},
\een
where we identify $\ze_{\mu}$ with the vector
$\ze^{}_{\Tc^{\tss0}}$ associated with the initial tableau
$\Tc^{\tss0}$ of the skew
shape $\Ga_{\la}/\mu$.
Hence, using the above formula for $g_{1,\mu}$ we find that
\ben
z_{m+1,1}\ts\ze^{}_{\Tc^{+}}=\frac{(-1)^{k-1}}{(r+m)^{k-1}}
{}\ts\prod_{p=1}^{k-1}(r+\la_{m+p}+m)
\ts\prod_{j=2}^m\ts
\Big(\frac{r-\si_j+1}
{r-\si_j}\Big)^{\min\{\la_j-\mu_j,\ts k-1\}}
\ts\ze_{\mu},
\een
where $\Tc^{+}$ is the tableau defined in \eqref{parta}.
We need to replace $\mu$ by $\mu-\de_1$ so that
$r$ is replaced by $r-1$ and $k=\la_1-r$ is replaced by
$k+1$. Assuming that $\Tc^{+}$ is now associated with
$\mu-\de_1$, the above relation takes the form
\ben
z_{m+1,1}\ts\ze^{}_{\Tc^{+}}=\frac{(-1)^{k}}{(r+m-1)^{k}}
\ts\prod_{p=1}^{k}(r+\la_{m+p}+m-1)
\ts\prod_{j=2}^m\ts
\Big(\frac{r-\si_j}
{r-\si_j-1}\Big)^{\min\{\la_j-\mu_j,\ts k\}}
\ts\ze_{\mu-\de_1}.
\een
To calculate $z_{m+1,1}\ts\ze^{}_{\Tc}$ for an arbitrary
skew tableau of shape $\Ga_{\la}/\mu$, we
use \eqref{bklact}
and the second relation of Lemma~\ref{lem:repbc}.
Applying appropriate operators
of the form $B_p(-l^{\tss\prime}_{r+p,j})$
to both sides of \eqref{zmponei}, we find that
\ben
f_{1,\Tc}=f_{1,\Tc^+}\ts
\prod_{p=1}^{n-1}\ts\prod_{j=1}^{r+p-1}\ts
\frac{r+l^{\tss0}_{r+p,j}+p-1}{r+l^{\tss\prime}_{r+p,j}+p-1},
\een
where the parameters $l^{\tss0}_{r+p,j}$ and
$l^{\tss\prime}_{r+p,j}$ are associated with
the tableaux $\Tc^{\tss0}$ and $\Tc$, respectively.
This completes the proof of \eqref{ekkmodd} and the theorem.
\epf

\bre\label{rem:nomalize}
Note that any normalization of the vectors $\ze_{\mu}$
with normalization constants depending only on $\mu$
(and $\la$) will only affect the formulas for the action
of $E_{m,m+1}$ and $E_{m+1,m}$ and leave the formulas for the
action of the
generators of the subalgebras
$\gl_m$ and $\gl_n$ of $\gl_{m|n}$ in the new basis
unchanged. This can be used to construct
basis vectors $\xi_{\La}=N(\La)\ts\ze_{\La}$ with simpler
expansions of $E_{m,m+1}\ts\xi_{\La}$ and
$E_{m+1,m}\ts\xi_{\La}$.
\qed
\ere

Theorem~\ref{thm:baslla} implies
the formula for the character of $L(\la)$ defined in \eqref{chardef},
which was originally found by Berele and Regev~\cite{br:hy}
and Sergeev~\cite{s:ta}.
Recall that given a Young diagram $\rho$,
the corresponding supersymmetric Schur polynomial
$s^{}_{\rho}(x)$ in the variables
$x=(x_1,\dots,x_m/x_{m+1},\dots,x_{m+n})$
is defined by the formula
\ben
s^{}_{\rho}(x)=
\sum_{\La}
\prod_{\al\in\rho}x^{}_{T(\al)},
\een
summed over the supertableaux $\La$ of shape $\rho$,
where $T(\al)$ denotes
the entry in the box $\al$ of the diagram $\rho$.
An alternative expression for $s^{}_{\rho}(x)$ is provided
by the Sergeev--Pragacz formula; see e.g. \cite[p.\ts61]{m:sf}.

\bco\label{cor:char}
The character $\ch L(\la)$ coincides with the supersymmetric
Schur polynomial
$s^{}_{\ts\Ga_{\la}}(x_1,\dots,x_m/x_{m+1},\dots,x_{m+n})$
associated with the Young diagram $\Ga_{\la}$.
\qed
\eco

\bex\label{ex:mone} In the case $n=1$, the basis
$\ze_{\La}$ of the representation
$L(\la_1,\dots,\la_m\ts |\ts \la_{m+1})$
of the Lie superalgebra $\gl_{m|1}$ can be parameterized
by the patterns
\begin{align}
&\qquad\la^{}_{m1}\qquad\la^{}_{m2}
\qquad\qquad\cdots\qquad\qquad\la^{}_{mm}\non\\
&\qquad\qquad\la^{}_{m-1,1}\qquad\ \ \cdots\ \
\ \ \qquad\la^{}_{m-1,m-1}\non\\
\Uc=&\quad\qquad\qquad\cdots\qquad\cdots\qquad\cdots\non\\
&\quad\qquad\qquad\qquad\la^{}_{21}\qquad\la^{}_{22}\non\\
&\quad\qquad\qquad\qquad\qquad\la^{}_{11}  \non
\end{align}
see Sec.~\ref{sec:mt}. The top row runs over partitions
$(\la^{}_{m1},\dots,\la^{}_{mm})$ such that either
$\la_{mj}=\la_j$ or $\la_{mj}=\la_j-1$ for each $j=1,\dots,m$.
The action of the generators
$E_{ss}$ with $s=1,\dots,m$ and the generators
$E_{s,s+1}$, $E_{s+1,s}$ with $s=1,\dots,m-1$
on the basis
vectors $\ze_{\ts\Uc}$ is given by the Gelfand--Tsetlin formulas
\begin{align}
E_{ss}\ts \ze^{}_{\ts\Uc}&=\Big(\sum_{i=1}^s\lambda^{}_{si}
-\sum_{i=1}^{s-1}\lambda^{}_{s-1,i}\Big)\ts
\ze^{}_{\ts\Uc},
\non
\\
E_{s,s+1}\ts \ze^{}_{\ts\Uc}
&=-\sum_{i=1}^s \frac{(l^{}_{si}-l^{}_{s+1,1})
\dots (l^{}_{si}-l^{}_{s+1,s+1})}
{(l^{}_{si}-l^{}_{s1})\ldots
\wedge_i\ldots(l^{}_{si}-l^{}_{ss})}
\ts
\ze^{}_{\ts\Uc+\de_{si}},
\non\\
E_{s+1,s}\ts \ze^{}_{\ts\Uc}
&=\sum_{i=1}^s \frac{(l^{}_{si}-l^{}_{s-1,1})
\dots (l^{}_{si}-l^{}_{s-1,s-1})}
{(l^{}_{si}-l^{}_{s1})\ldots
\wedge_i\ldots(l^{}_{si}-l^{}_{ss})}
\ts
\ze^{}_{\ts\Uc-\de_{si}},
\non
\end{align}
where the array $\Uc\pm\de_{si}$ is obtained from $\Uc$
by replacing $\la_{si}$ with $\la_{si}\pm 1$, and
$\ze^{}_{\ts\Uc}$ is considered to be equal to zero if
the array $\Uc$ is not a pattern. For the action
of the generators $E_{m+1,m+1}$, $E_{m,m+1}$ and $E_{m+1,m}$
we have
\begin{align}
E_{m+1,m+1}\ts \ze^{}_{\ts\Uc}&=\Big(\sum_{i=1}^{m+1}\lambda^{}_{i}
-\sum_{i=1}^{m}\lambda^{}_{mi}\Big)\ts
\ze^{}_{\ts\Uc},
\non\\[1em]
E_{m,m+1}\ts \ze^{}_{\ts\Uc}
&=\sum_{i=1}^m
(l_{mi}+\la_{m+1}+m)
\non\\
{}&\times{}\ts\prod_{j=1}^{i-1}
(-1)^{\la_j-\la_{mj}}\ts\frac{l_{mi}-l_j}{l_{mi}-l_{mj}}\ts
\prod_{j=i+1,\ \la_j-\la_{mj}=1}^{m}
\frac{l_{mi}-l_{mj}+1}{l_{mi}-l_j+1}
\ts
\ze^{}_{\ts\Uc+\de_{mi}},
\non\\[1em]
E_{m+1,m}\ts \ze^{}_{\ts\Uc}
&=\sum_{i=1}^m
\frac{(l^{}_{mi}-l^{}_{m-1,1})
\dots (l^{}_{mi}-l^{}_{m-1,\tss m-1})}
{(l^{}_{mi}-l^{}_{m1})\ldots\wedge_i\ldots(l^{}_{mi}-l^{}_{mm})}
\non\\
{}&\times{}\prod_{j=1}^{i-1}
(-1)^{\la_j-\la_{mj}}\ts\frac{l_{mi}-l_{mj}-1}{l_{mi}-l_j-1}
\ts\prod_{j=1,\ \la_{j}-\la_{mj}=1}^{i-1}
\frac{l_{mi}-l_{mj}}{l_{mi}-l_j}\ts
\ze^{}_{\ts\Uc-\de_{mi}},
\non
\end{align}
where the coefficients in the expansion of
$E_{m,m+1}\ts \ze^{}_{\ts\Uc}$ were simplified
with the use of Lemma~\ref{lem:ida}.
The basis $\ze^{}_{\ts\Uc}$ of $L(\la_1,\dots,\la_m\ts |\ts \la_{m+1})$
coincides with that of \cite{p:if}
up to a normalization.
\qed
\eex

\bex\label{ex:onen} In the case $m=1$, the basis
$\ze_{\La}$ of the representation
$L(\la_1\ts |\ts \la_2,\dots,\la_{n+1})$
of the Lie superalgebra $\gl_{1|n}$ can be parameterized
by the trapezium patterns
\begin{align}
\la^{\prime}_{r+n,1}\quad\la^{\prime}_{r+n,2}\qquad\cdots
\qquad\qquad&\cdots\qquad\qquad\qquad\la^{\prime}_{r+n,r+n}\non\\
\Vc=\qquad\qquad\dddots\quad\dddots
\qquad\cdots\qquad\quad&\cdots\qquad\qquad\qquad\antiddots\non\\
\qquad\qquad\la^{\prime}_{r+1,1}\quad \la^{\prime}_{r+1,2}
\qquad\quad&\cdots\quad\qquad\la^{\prime}_{r+1,r+1}\non\\
\quad\qquad\qquad\qquad 1
\qquad\ \ 1
\qquad\quad&\cdots\qquad 1\non
\end{align}
where the number $r$ of $1$'s in the bottom row
is nonnegative and varies
between $\la_1-n$ and $\la_1$. The top row coincides with
$(\la'_1,\dots,\la'_q,0,\dots,0)$, where $q=\la_1$.
The action of the generators
$E_{ss}$ with $s=1,\dots,n+1$ and the generators
$E_{s,s+1}$, $E_{s+1,s}$ with $s=2,\dots,n$
on the basis
vectors $\ze_{\ts\Vc}$ is given by the Gelfand--Tsetlin formulas
\eqref{ekk}--\eqref{ekkm}. The formulas for the action
of $E_{12}$ and $E_{21}$ are given in \eqref{ekkpodd} and
\eqref{ekkmodd}, respectively.

Specializing further and taking $n=2$ with $\la_1\geqslant 2$
we can parameterize the basis vectors $\ze_{\Vc}$
of the $\gl_{1|2}$-module
$L(\la_1\ts|\ts\la_2,\la_3)$ by trapezium patterns $\Vc$
of four types:
\begin{align}
&\ \la_2+1\quad\la_3+1\quad1
\quad\cdots\quad1\quad0\quad0\non\\
\Vc_a^{(1)}\quad=\quad&\ \ \quad\quad a+1\quad\quad 1
\quad\cdots\quad1\quad1\quad0\non\\
&\quad\quad\quad\quad\quad\quad1
\quad\cdots\quad1\quad 1\quad 1\non
\end{align}
the number of $1$'s in the bottom row is $\la_1$;
\begin{align}
&\ \la_2+1\quad\la_3+1\quad1
\quad\cdots\quad1\quad0\non\\
\Vc_a^{(2)}\quad=\quad&\ \ \quad\quad a+1\quad\quad 1
\quad\cdots\quad1\quad1\non\\
&\quad\quad\quad\quad\quad\quad1
\quad\cdots\quad1\quad 1\non
\end{align}
the number of $1$'s in the bottom row is $\la_1-1$;
\begin{align}
&\ \la_2+1\quad\la_3+1\quad1
\quad\cdots\quad1\quad0\non\\
\Vc_a^{(3)}\quad=\quad&\ \ \quad\quad a+1\quad\quad 1
\quad\cdots\quad1\quad0\non\\
&\quad\quad\quad\quad\quad\quad1
\quad\cdots\quad1\quad 1\non
\end{align}
the number of $1$'s in the bottom row is $\la_1-1$;
\begin{align}
&\ \la_2+1\quad\la_3+1\quad1
\quad\cdots\quad1\non\\
\Vc_a^{(4)}\quad=\quad&\ \ \quad\quad a+1\quad\quad 1
\quad\cdots\quad1\non\\
&\quad\quad\quad\quad\quad\quad1
\quad\cdots\quad1\non
\end{align}
the number of $1$'s in the bottom row is $\la_1-2$, where the
parameter $a$ in all cases runs over the integers such that
$\la_3\leqslant a\leqslant \la_2$. Thus,
$\dim L(\la_1\ts|\ts\la_2,\la_3)=4\tss(\la_2-\la_3+1)$.
The formulas for the action of the odd
generators of the Lie superalgebra $\gl_{1|2}$
provided by Theorem~\ref{thm:baslla} have the form
\ben
\bal
E_{12}\ts\ze_{\Vc_a^{(1)}}&=0,\qquad
E_{12}\ts\ze_{\Vc_a^{(2)}}
=\frac{(a+\la_1)(a+\la_1+1)}{\la_1+\la_2+1}\ts\ze_{\Vc_a^{(1)}}\\
E_{12}\ts\ze_{\Vc_a^{(3)}}&=0,\qquad
E_{12}\ts\ze_{\Vc_a^{(4)}}=
-\frac{(\la_1-1)(a+\la_1)(a+\la_1-1)}{(\la_1+\la_2)(\la_1+\la_2-1)}
\ts\ze_{\Vc_a^{(3)}}
\eal
\een
and
\ben
\bal
E_{21}\ts\ze_{\Vc_a^{(2)}}&=0,\qquad
E_{21}\ts\ze_{\Vc_a^{(1)}}
=\frac{\la_1+\la_2+1}{a+\la_1+1}\ts\ze_{\Vc_a^{(2)}}\\
E_{21}\ts\ze_{\Vc_a^{(4)}}&=0,\qquad
E_{21}\ts\ze_{\Vc_a^{(3)}}=
-\frac{(\la_1+\la_2)(\la_1+\la_2-1)}{(\la_1-1)(a+\la_1)}
\ts\ze_{\Vc_a^{(4)}}.
\eal
\een
\qed
\eex

\end{document}